\documentclass[a4paper]{article}

\usepackage{amsmath,amssymb,amsthm}
\usepackage{mathrsfs}

\usepackage[matrix,arrow,curve]{xy}
\usepackage[left=2cm,right=2cm,top=2cm,bottom=2cm,bindingoffset=0cm]{geometry}
\usepackage{graphicx}
\usepackage{nccfoots}
\usepackage[hang,flushmargin]{footmisc} 
\usepackage{hyperref}

\usepackage{pgffor}

\usepackage{pgf,tikz}
\usepackage{mathrsfs}
\usetikzlibrary{arrows,calc,intersections,through,shapes}

\newtheorem{remark}{Remark}
\newtheorem{prop}{Proposition}
\newtheorem{lemma}{Lemma}
\newtheorem{cor}{Corollary}

\newtheorem{observ}{Observation}
\newtheorem{theo}{Theorem}
\newtheorem{problem}{Problem}
\newtheorem{definition}{Definition}
\newtheorem{conjecture}{Conjecture}

\DeclareMathOperator{\Int}{Int}

\DeclareMathOperator{\iopta}{\iota}
\newcommand{\HH}{{\mathcal H}}

\renewcommand{\H}{\HH^1}

\newcommand{\St}{{\mathcal St}}

\renewcommand{\H}{\HH^1}

\newcommand{\forget}[1]{}

\def\eps{\varepsilon}

\def\Rr{\mathcal{R}}
\def\Rg{\mathcal{R}_{\mathrm{geoemb}}}
\def\Rc{\mathcal{R}_{\mathrm{lencrit}}}

\def\Aa{\mathcal{A}}
\def\Rm{\mathcal{R}_{\mathrm{locmin}}}

\def\PP{\mathbb{P}}

\def\Oo{\mathcal{O}}

\title{On uniqueness in Steiner problem}

\author{M. Basok$^{2}$, D. Cherkashin$^{1,3}$, N. Rastegaev$^1$, Y. Teplitskaya$^{1,4}$}

\begin{document}

\maketitle

\begin{abstract}
We prove that the set of $n$-point configurations for which the solution of the planar Steiner problem is not unique has the Hausdorff dimension at most $2n-1$ (as a subset of $\mathbb{R}^{2n}$). Moreover, we show that the Hausdorff dimension of the set of $n$-point configurations on which at least two locally minimal trees have the same length is also at most $2n-1$.
Methods we use essentially require rely upon the theory of subanalytic sets developed in~\cite{bierstone1988semianalytic}. Motivated by this approach we develop a general setup for the similar problem of uniqueness of the Steiner tree where the Euclidean plane is replace by an arbitrary analytic Riemannian manifold $M$. In this setup we argue that the set of configurations possessing two locally-minimal trees of the same length either has the dimension $n\dim M-1$ or has a non-empty interior. We provide an example of a two-dimensional surface for which the last alternative holds.

In addition to abovementioned results, we study the set of set of $n$-point configurations for which there is a unique solution of the Steiner problem in $\mathbb{R}^d$. We show that this set is path-connected.
\end{abstract}

\section{Introduction}
\Footnote{}{\noindent $^1$ Chebyshev Laboratory, St. Petersburg State University, 14th Line V.O., 29B, Saint Petersburg 199178 Russia\\
$^2$ University of Helsinki, Finland\\ 
$^3$ \noindent Institute of Mathematics and Informatics, Bulgarian Academy of Sciences\\ 
$^4$ \noindent Mathematical Institute, Leiden University, the Netherlands}

We consider the following form of the Steiner tree problem:
\begin{problem}\label{Problem1}
For a given finite set $P = \{x_1,\dots ,x_n\} \subset \mathbb{R}^d$ find a connected set $\St$ with minimal length (one-dimensional Hausdorff measure $\H$) containing $P$.
\end{problem}

Throughout the article $n \geq 4, d \geq 2$ are natural numbers. All the solutions of Problem~\ref{Problem1} for $n \leq 3$ are known in the explicit form 
since 17-th century.

A solution of Problem~\ref{Problem1} is called \textit{Steiner tree}.
It is known that such an $\St$ always exists (but is not necessarily unique, see Fig.~\ref{fig1}) and that it is a union of a finite set of segments. Thus, $\St$ can be represented as a graph, embedded into the Euclidean space, such that its set of vertices contains $P$ and all its edges are straight line segments. This graph is connected and does not contain cycles, i.e. is a tree, which explains the naming of $\St$. It is known that the maximal degree of the vertices of $\St$ is at most $3$. Moreover, only vertices $x_i$ can have degree $1$ or $2$, all the other vertices have degree $3$ and are called \textit{Steiner points} while the vertices $x_i$ are called \textit{terminals}. 
Vertices of the degree $3$ are called \textit{branching points}.
The angle between any two adjacent edges of $\St$ is at least $2\pi/3$.
That means that for a branching point the angle between any two segments incident to it is exactly $2\pi/3$, and these three segments belong to the same 2-dimensional plane.

The number of Steiner points in $\St$ does not exceed $n-2$. A Steiner tree with exactly $2n-2$ vertices is called \textit{full}.  
Every terminal point of a full Steiner tree has degree one.

For a given finite set $P \subset \mathbb{R}^d$ consider a connected acyclic set $S$ containing $P$.
Then $S$ is called a \textit{locally minimal tree} if $\overline{S \cap B_\varepsilon (x)}$ is a Steiner tree for 
$(\{x\} \cap P) \cup \partial(S \cap B_\varepsilon (x))$ for every point $x \in S$ and small enough $\varepsilon>0$. Clearly every Steiner tree is locally minimal and not vice versa.
Locally minimal trees have all the mentioned properties of Steiner trees except the minimal length condition. 
So locally minimal trees inherit the definitions of terminals, Steiner points, branching points and fullness.
Proof of the listed properties of Steiner and locally minimal trees together with an additional information on them can be found in
book~\cite{hwang1992steiner} and in article~\cite{gilbert1968steiner}.

Similar problems could also be considered in abstract metric spaces. In the most general form the problem would be to connect a set (not necessarily finite or countable) of subsets of an arbitrary metric space in a minimal way with respect to the metric~\cite{paolini2013existence}, see~Section~\ref{subsec:arbitrary_metric_space}.

The Steiner problem may have several solutions starting with $n = 4$ (see Fig.~\ref{fig1}). 
The main result of the paper implies the uniqueness of a solution for a general input. 

Let us denote by $\mathbb{P}_d := (\mathbb{R}^d)^n \setminus diag$ the space of labeled $n$-point configurations $x_1, \dots , x_n \in \mathbb{R}^d$ of distinct points in the Euclidean space, where $diag$ is the union of $(dn - d)$-dimensional subspaces $x_i = x_j$, $i \neq j$.
Note that every point of $\mathbb{P}_d$ corresponds to some labeled non-degenerate configuration; so let us consider $\mathbb{P}_d$ as a configuration space.

%add local ambiguous
A configuration $P \in \mathbb{P}_d$ is \textit{ambiguous} if there are several Steiner trees for $P$. 
Ivanov and Tuzhilin proved~\cite{ivanov2006uniqueness} that the complement to the set of ambiguous configurations contains an open dense subset\Footnote{1}{For some reason they call it ``general position''} of $\mathbb{P}_2$. 
Edelsbrunner and Strelkova~\cite{edelsbrunner2012configuration} asked whether the measure of ambiguous configurations is zero or not. 
We provide a positive answer by proving the following stronger statement.

\begin{theo}
Assume that $n\geq 4$. Then the set of planar ambiguous configurations in $\mathbb{P}_2$ has the Hausdorff dimension $2n - 1$.
\label{main}
\end{theo}

\begin{figure}[h]
    \centering
    \begin{tikzpicture}
    \def\r{1.5cm}
    \draw[ultra thick, blue]
        (-\r, \r) coordinate(x1) node[black, above right]{$1$} --++ (-60:{\r/cos(30)}) coordinate (x5);
    \draw[ultra thick, blue]
        (\r,\r) coordinate(x2) node[black, above left]{$2$} --++ (-120:{\r/cos(30)}) coordinate (x6);
    \draw[ultra thick, blue]
        (\r, -\r) coordinate(x3) node[black, below left]{$3$} --++ (120:{\r/cos(30)});
    \draw[ultra thick, blue]
        (-\r,-\r) coordinate(x4) node[black, below right]{$4$} --++ (60:{\r/cos(30)});
    \draw[ultra thick, blue]
        (x5) -- (x6);
    \foreach \x in{1,2,...,6}{
        \fill (x\x) circle (2pt);
    }
\end{tikzpicture}
\hspace{2cm}
\begin{tikzpicture}
    \def\r{1.5cm}
    \draw[ultra thick, blue]
        (-\r, \r) coordinate(x1) node[black, above right]{$1$} --++ (-30:{\r/cos(30)}) coordinate (x5);
    \draw[ultra thick, blue]
        (\r,\r) coordinate(x2) node[black, above left]{$2$} --++ (-150:{\r/cos(30)});
    \draw[ultra thick, blue]
        (\r, -\r) coordinate(x3) node[black, below left]{$3$} --++ (150:{\r/cos(30)}) coordinate (x6);
    \draw[ultra thick, blue]
        (-\r,-\r) coordinate(x4) node[black, below right]{$4$} --++ (30:{\r/cos(30)});
    \draw[ultra thick, blue]
        (x5) -- (x6);
    \foreach \x in{1,2,...,6}{
        \fill (x\x) circle (2pt);
    }
\end{tikzpicture}
    \caption{An example of nonunique solution. Labelled points form a square.}
    \label{fig1}
\end{figure}

\subsection{Topology and embedding class of a tree}
\label{subsec:topology_def}

For the sake of convenience and completeness, we would like to begin our discussion with a careful introduction of the concept of ``topology'' of a tree ofter used in the context of the Steiner problem. As it is usually done in the literature (see, for instance~\cite{gilbert1968steiner,hwang1992steiner}), we define the \textit{topology} of a Steiner tree to be the corresponding abstract topological graph with labelled terminals and unlabelled Steiner points. Thus, a topology $T$ is a topological space with a tree structure, and some vertices of $T$, including all its leaves and vertices of degree 2, are labelled. Moreover, we assume that all vertices of $T$ have degrees at most 3, as it naturally holds for any Steiner tree.

Note that two trees embedded in a different (non-homotopic) way into the plane may have the same topology.
To distinguish non-homotopic embeddings Edelsbrunner and Strelkova~\cite{edelsbrunner2019configuration} 
introduced another invariant way to describe the topological type of the tree which we call the ``embedding class''.
Below we introduce several ways to define the embedding class of a tree commonly used. We include the proof of their equivalence in the Appendix for the sake of completeness.

Let $T$ be a combinatorial tree. Let $\vec{E}(T)$ denote the set of oriented edges of $T$ (in particular, $|\vec{E}(T)| = 2|E(T)|$). Given an edge $\vec{e}\in \vec{E}$ denote by $o(\vec{e})$ the origin and by $t(\vec{e})$ the tail. Let us say that a bijection $\sigma: \vec{E}(T)\to \vec{E}(T)$ determines \emph{a cyclic order around each vertex of $T$} if $o(\sigma(\vec{e})) = o(\vec{e})$ for any $\vec{e}$ and for any vertex $v\in V(T)$ and $\vec{e}$ such that $o(\vec{e}) = v$ the set $\vec{e}, \sigma(\vec{e}),\sigma^2(\vec{e}),\dots$ is exactly the set of oriented edges emanating from $v$.

The following classical lemma defines the embedding class:

\begin{lemma}
    \label{lemma:class_of_embeddings}
    Let a positive integer $n$ be fixed. The following three sets are in natural bijection:
    \begin{enumerate}
        \item The set $PM_1$ of pairs $(T,\sigma)$, where $T$ is a combinatorial tree with all vertices of degree at most 3, with $n$ labelled vertices, including all leaves of $T$ and vertices of degree 2, and $\sigma:\vec{E}(T)\to \vec{E}(T)$ is a bijection determining a cyclic order around each vertex.
        \item The set $PM_2$ of pairs $(T,[f])$, where $T$ is a combinatorial tree with all vertices of degree at most 3, with $n$ labelled vertices, including all leaves of $T$ and vertices of degree 2, $f$ is a bijection between $\vec{E}(T)$ and the set $\partial D$ of edges of the regular $|\vec{E}(T)|$-gon $D$ (0-gon is assumed to be empty) oriented clockwise such that if $o(f(\vec{e}_1))=t(f(\vec{e}_2))$, then $o(\vec{e}_1) = t(\vec{e}_2)$, and $[f]$ is the equivalence class of $f$ with respect to the cyclic shift on $\partial D$. 
        \item The set $PM_3$ of pairs $(T,[\iopta])$, where $T$ is a topology with $n$ labelled vertices, $\iopta$ is some embedding of $T$ into the plane and $[\iopta]$ is the homotopy class of $\iopta$ in the space of embeddings.
    \end{enumerate}
\end{lemma}
\begin{remark}
    Each of these three sets can be considered as the set of \emph{plane maps} or \emph{ribbon graphs} with the tree-like skeletons and some labelled vertices (see~\cite{lando2010graphs}). 
\end{remark}

We prove Lemma~\ref{lemma:class_of_embeddings} in the Appendix.

Note that the regular polygon $D$ from the set $PM_2$ naturally corresponds to the outer face of the planar graph $\iopta(T)$ for $\iopta$ coming from $PM_3$. Using Lemma~\ref{lemma:class_of_embeddings}, we identify $PM_1,PM_2$ and $PM_3$, so that given, say $(T,\sigma)\in PM_1$ we will always assume that we are also given the corresponding $(T,[f])\in PM_2$ and $(T,[\iopta])\in PM_3$ and will use the corresponding notation if it does not lead to a confusion.

Now, we introduce another (fourth) way to encode embeddings of a topological tree, which was originally used by Edelsbrunner and Strelkova.

Let $(T,\sigma)\in PM_1$ and $\varphi = \alpha\circ \sigma$, where  $\alpha: \vec{E}(T)\to \vec{E}(T)$ is the involution reversing the orientation. Assume that $T$ has $n$ labels. Let $A = \{1,2,3,\dots,n,b\}$ be the alphabet on $n+1$ letters, $n$ of them are numbers from $1$ to $n$, and $(n+1)$-th is the special letter $b$. Let
\[
    \mathbf{C} = \left(\bigcup_{k\geq 0}A^k\right)/_\text{cyclic shift}
\]
be the set of all words build from this alphabet considered up to the cyclic shift. Let $\vec{e}\in \vec{E}(T)$ be arbitrary and $N = |\vec{E}(T)|$. Then, given $(T,\sigma)$, define the word $C(T,\sigma)\in \mathbf{C}$ by the following rule: fix a vector $\vec{e}_0\in \vec{E}(T)$ and set $\vec{e}_i = \varphi(\vec{e}_{i-1})$, then define $C(T,\sigma) = a_0a_1\dots a_{N-1}$, where $a_i$ is the label of $o(\varphi(\vec{e}_i))$ if $o(\varphi(\vec{e}_i))$ is labelled, and $a_i = b$ else; if $T$ consists of one vertex, then the word $C(T,\sigma)$ is the empty word.

Let $(T,[f])\in PM_2$ corresponds to $(T,\sigma)$ and $D$ be the regular $N$-gon. As $D$ can be seen as the outer face of the planar graph $\iopta(T)$, there is a many-to-one correspondence between the vertices of $D$ and the vertices of $T$. Then the word $C(T,\sigma)$ is nothing but the list of vertices obtained by going along the boundary of $D$; each time we met a vertex walking along $\partial D$, we add its label to $C(T,\sigma)$, or the letter $b$ if the vertex does not have a label. For example, we have $C = 1bb2b3bb4b$ and $C = 1b2bb3b4bb$ (and we could also write $C = b2bb3b4bb1$ in the latter case as we factorized by a cyclic shift) for the left and the right trees on the Fig.~\ref{fig1} respectively.

\begin{lemma}
    \label{lemma:def_of_embedding_type}
    The morphism $(T,\sigma)\mapsto C(T,\sigma)$ is injective from the set of pairs $(T,\sigma)$ to $\mathbf{C}$.
\end{lemma}

We prove Lemma~\ref{lemma:def_of_embedding_type} in the Appendix.

Starting from now we will call the word $C(T,\sigma)$ an \emph{embedding class}. Using Lemma~\ref{lemma:class_of_embeddings} and Lemma~\ref{lemma:def_of_embedding_type} we will feel free to identify the embedding class with the homotopy class of embeddings defined in several ways presented in aforementioned lemmas.

\subsection{Connectedness in \texorpdfstring{$\mathbb{P}_d$}{Pd}}

%A \textit{topology} $T$ of a labeled Steiner or locally minimal tree $S$ is the \textcolor{red}{clockwise cyclic} order of the vertices of $S$ (including Steiner points, whose all are denoted by $b$). For example Steiner trees for a square have topologies $1bb2b3bb4b$ and $1b2bb3b4bb$, see the left and right parts of Fig.~\ref{fig1}, respectively.
Let us return to our analysis of Steiner trees. We say that a topology $T$ of a tree $S$ is \textit{full} if the corresponding tree is full. %(obviously it does not depend on a tree).
%The \textit{degree} of a label $A$ in a topology $T$ is the number of times $A$ appears in the cyclic order.
%Clearly, the degree of a label in $T$ coincides with the graph-theoretical degree of the corresponding vertex of $S$.
Further, let us call a topology $T$ \textit{realizable} for a set $P$ if there exists such a locally minimal tree $S(P)$ with topology $T$; we will denote this tree by $S_T(P)$.

\begin{prop}[Melzak,~\cite{melzak1961problem}]
If a topology $T$ is realizable for $P \in \mathbb{P}$ then the realization $S_T(P)$ is unique.
\label{melzakuniq}
\end{prop}
Proposition~\ref{melzakuniq} shows that $S_T(P)$ is uniquely defined.
Moreover one can construct (or show that it is impossible) $S_T(P)$ in a linear time~\cite{hwang1986linear}. 
However, we rarely know a priori, which topology gives a Steiner tree. 
Although the number of possible topologies for an $n$ points configuration is finite, checking all of them may consume a lot of time, since this number of topologies grows very fast with $n$, see~\cite{gilbert1968steiner,hwang1992steiner}. Indeed, the Steiner tree problem is NP-complete~\cite{garey1977complexity}.

We need the following generalization of Proposition~\ref{melzakuniq}. For a full topology $T$ define $D(T)$ as the set
of topologies that can be obtained from $T$ by shrinking some edges connecting a terminal with a Steiner point (these edges should have pairwise different ends).

\begin{prop}[Gilbert--Pollak~\cite{gilbert1968steiner}, Hwang--Weng~\cite{hwang1992shortest}]
Let $T$ be a full topology and $P \in \mathbb{P}$. Consider function $L(y_1,\dots , y_{n-2}) : (\mathbb{R}^2)^n \to \mathbb R$ which is the length of a tree on the vertex set $P \cup \{y_1,\dots, y_{n-2}\}$ with straight edges and topology $T$ (we allow $y_i$ coincide with terminals). Then $L$ has unique local minimum and so there is exactly one realization with a topology from $D(T)$.
\label{GPuniq}
\end{prop}

A \textit{generic topology} is a topology without terminals of degree 3.

\begin{observ}
\begin{itemize}
    \item[(i)] Every generic topology $R$ belongs to exactly one set $D(T)$, because the reverse procedure (replacing every vertex $A$ of degree 2 in $R$ on a Steiner point $b$ and add edge $bA$) leads to a full topology $T$.

    \item[(ii)] Suppose that $\St$ is the unique Steiner tree for some $P \in \mathbb{P}$ and has a generic topology $R \in D(T)$ for some full topology $T$. Then for some positive $\eta > 0$ and any other full topology $T'$ the length of the realization from $D(T')$ exceeds $\HH(\St)$ by at least $\eta$.
    If one changes every terminal with a point from its $\eta/n$-neighborhood, then by triangle inequality
    a perturbed configuration $P'$ has a unique Steiner tree $\St(P')$ and the topology of $\St(P')$ belongs to $D(T)$.
    
    \item[(iii)] Configurations $P \in \mathbb{P}_2$ for which there is a locally minimal tree with non-generic topology has the Hausdorff dimension $2n-2$.
    
\end{itemize}
\label{obs:topologies}
\end{observ}

In the paper we study the way realizations and minimal realizations of different embedding classes divide the configuration space. 

A similar research topic appears in~\cite{edelsbrunner2012configuration,edelsbrunner2019configuration}, where the connectedness of some sets related to an embedding class $EC$, is studied.
Let $\Omega(EC)$ be a subset of $\mathbb{P}$ consisting of all $P \in \mathbb{P}$ for which $EC$ is realizable.
Note that for every embedding class $EC$ the set $\Omega(EC)$ is path-connected.

\begin{theo}[Edelsbrunner--Strelkova,~\cite{edelsbrunner2012configuration,edelsbrunner2019configuration}]
Let $EC$ be an embedding class. Then the subset of $\mathbb{P}_d$ for which a Steiner tree is unique and has the embedding class $EC$ is path-connected.
\label{uniquepathconnectedness}
\end{theo}

In the planar case they also obtained the following result. 

\begin{theo}[Edelsbrunner--Strelkova,~\cite{edelsbrunner2012configuration,edelsbrunner2019configuration}]
Let $EC$ be a full embedding class. Then the subset of $\mathbb{P}_2$ for which a Steiner tree has the embedding class $EC$ is path-connected.
\end{theo}

The second result of our paper is the following.

\begin{theo}
The subset of $\mathbb{P}$ for which there is a unique Steiner tree is path-connected.
\label{second}
\end{theo}
The proof of Theorem~\ref{second} is constructive (modulo Theorem~\ref{uniquepathconnectedness}) and the embedding classes of the Steiner tree is known at every point of a constructed path.

\textbf{Structure of the paper.} Section~\ref{can} devoted to the proof of Theorem~\ref{second}. We begin with Section~\ref{subsec:universal_ST} where we recall the construction of a universal Steiner tree introduces by Paolini, Stepanov and Teplitskaya. Using this construction we introduce a canonical realization of an embeddin class in Section~\ref{subsec:can_real} and use it to prove Theorem~\ref{second} in Section~\ref{subsec:proof_of_second}.

We proceed with Section~\ref{sec:proof of main theorem}. In Section~\ref{subsec:Subanalytic subsets of a real analytic manifold} we review some facts about subanalytic subsets following~\cite{bierstone1988semianalytic}. We prove Theorem~\ref{main} in Section~\ref{subsec:proof_of_Theorem_on_dimension}.

Finally, in Section~\ref{sec:Steiner trees in analytic Riemannian manifolds} is briefly discuss Steiner problem on an arbitrary analytic manifold. We introduce several configuration spaces and study their local structures in Sections~\ref{subsec:arbitrary_metric_space}--\ref{subsec:Rlocmin}. In Section~\ref{subsec:topological_example} we provide an example of a surface where the set of configurations having two locally minimal trees of the same length has a non-empty interior.

\section{Connectivity of the subset of \texorpdfstring{$\mathbb P$}{P} with a unique Steiner tree}
\label{can}

\subsection{A universal Steiner tree}
\label{subsec:universal_ST}

In this subsection we provide the construction of a unique Steiner tree with an infinite number of Steiner points from~\cite{paolini2015example} (a partial improvement of the result appeared in~\cite{cherkashin2023self}). 

Let $S_\infty$ be an infinite tree with vertices $y_0, y_1, y_2,\dots$ and edges given by $y_0y_1$ and $y_ky_{2k},\ y_ky_{2k+1},\ k\geq 1$. Thus, $S_\infty$ is an infinite binary tree with an additional vertex $y_0$ attached to the common parent $y_1$ of all other vertices $y_k,k\geq 2$. The goal of~\cite{paolini2015example} is to embed $S_\infty$ in the plane in such a way that the image of each finite subtree of $S_\infty$ will be the unique Steiner tree for the set of its vertices having degree 1 or 2. We define the embedding below by specifying the positions of $y_0,y_1,y_2,\dots$ on the plane.

Let $\Lambda = \{\lambda_i\}_{i=0}^\infty$ be a sequence of positive real numbers.
Define an embedding $\Sigma (\Lambda)$ of $S_\infty$ as a rooted binary tree with the root $y_0 = (0,0)$ the first descendant $y_1 = (1,0)$ and the ratio between edges of $(i+1)$-th and $i$-th levels being $\lambda_i$.
For a small enough $\{\lambda_i\}$ the set $\Sigma (\Lambda)$ see Fig.~\ref{pic:univerasaltree}.

\begin{figure}[h]
\begin{center}
\begin{tikzpicture}
  % points
  \def\lambda{0.3}
  \def\ty{0.4} %% 1-2*\lambda
  \def\r{7cm}
  \def\rr{\lambda*\r}
  \def\rrr{\ty*\ty*\r}

    \path (-2*\r,0) coordinate [label=above left:$y_0$] (y0);
    \path (-\r,0) coordinate [label=above left:$y_1$] (y1);

	\path (y1) ++(60:\ty*\r) coordinate [label=left:$y_2$] (y2);
	\path (y1) ++(-60:\ty*\r) coordinate [label=left:$y_3$] (y3);
	\path (y2) ++(120:\ty*\rr) coordinate [label=right:$y_4$] (y4);
	\path (y2) ++(0:\ty*\rr) coordinate [label={[xshift=-0.1cm, yshift=0.1cm]:$y_5$}] (y5);
	\path (y3) ++(0:\ty*\rr) coordinate [label={[xshift=-0.2cm, yshift=-0.42cm]:$y_6$}] (y6);
	\path (y3) ++(-120:\ty*\rr) coordinate [label=right:$y_7$] (y7);
	\path (y4) ++(180:\ty*\rrr) coordinate (y8);
	\path (y4) ++(60:\ty*\rrr) coordinate (y9);
	\path (y5) ++(60:\ty*\rrr) coordinate (y10);
	\path (y5) ++(-60:\ty*\rrr) coordinate (y11);
	\path (y6) ++(60:\ty*\rrr) coordinate (y12);
	\path (y6) ++(-60:\ty*\rrr) coordinate (y13);
	\path (y7) ++(-60:\ty*\rrr) coordinate (y14);
	\path (y7) ++(180:\ty*\rrr) coordinate (y15);

  % fork:
  \def\fork{\path[draw, color=blue, ultra thick]}

    \fork (y2) -- (y1) -- (y3);
	\fork (y4) -- (y2) -- (y5);
	\fork (y6) -- (y3) -- (y7);

   \path[draw, color=blue, ultra thick] (y0) -- (y1);

  \foreach \n [evaluate=\n as \m using {int(2*\n)}] in {4,5,6,7} { 
    \path[draw, color=blue, ultra thick, dotted] (y\n) -- (y\m);
  }

 \foreach \n [evaluate=\n as \m using {int(2*\n+1)}] in {4,5,6,7} { 
    \path[draw, color=blue, ultra thick, dotted] (y\n) -- (y\m);
  }
 
  \foreach \n in {0,1,2,3,4,5,6,7} {
    \fill [black] (y\n) circle (2pt);
  }

\end{tikzpicture}
\end{center}
\caption{Three iterations in the construction of $\Sigma_\infty$. The set $\Sigma_3$ is thick blue.}
\label{pic:univerasaltree}
\end{figure}

Let $A_\infty (\Lambda)$ be the union of the set of all leaves (limit points) of $\Sigma (\Lambda)$ and $\{y_0\}$.

\begin{theo}[Paolini--Stepanov--Teplitskaya,~\cite{paolini2015example}]
\label{th_steinMain0}
A binary tree $\Sigma(\Lambda)$ is a unique Steiner tree for $A_\infty (\Lambda)$ provided by $\lambda_i < 1/5000$ and $\sum_{i=1}^\infty \lambda_i  <\pi/5040$.
\end{theo}

We will use the following corollary of Theorem~\ref{th_steinMain0}, which explains why a full binary Steiner tree is universal, i.e. it contains a subtree with a given combinatorial structure.

\begin{cor}
In the conditions of Theorem~\ref{th_steinMain0} each connected closed subset $S$ of $\Sigma_\infty$ contained in $\Sigma_k$ for some $k$ has a natural tree structure. Moreover, every such an $S$ is the unique Steiner tree for any set $P$ containing the set of the vertices with the degree $1$ and $2$ of $S$.
\label{coruniquesubtree}
\end{cor}
\begin{proof}
    Let $S\subset \Sigma_\infty$ and $P\subset S$ satisfy the conditions of the corollary. The fact that $S$ is a tree is straightforward. Let $S'\neq S$ be any Steiner tree for $S$ and assume that $\H(S') \leq \H(S)$. Then it is clear that $\H((\Sigma_k\smallsetminus S)\cup S') \leq \H(\Sigma_k)$, but on the other hand $\{y_0\}\cup A_k\subset (\Sigma_k\smallsetminus S)\cup S'$, which contradicts to Theorem~\ref{th_steinMain0}.
\end{proof}

\subsection{Canonical realization of an embedding class}
\label{subsec:can_real}

Using the construction of the tree $\Sigma_\infty$ from the previous section we define the \emph{canonical realization} tree $\St_{EC}$ for any embedding class $EC$.

Fix a topological tree $S$ with the embedding class $EC$ and pick some vertex $v$ of $S$ of degree one. Then identify $S$ with the subtree of $\Sigma_\infty$ by mapping $v$ to the root $y_0$ of $\Sigma_\infty$ and mapping all the other vertices following the steps of the breadth-first search algorithm started from $v$, where at every vertex of degree 2 of $S$ we choose the left direction in $\Sigma_\infty$ (i.e. map the only child to $y_{2k}$ if the parent was mapped to $y_k$).
%By the construction every vertex of $\St_T$ coincides with $y_j$ for some $j$.

\subsection{Proof of Theorem~\ref{second}}
\label{subsec:proof_of_second}

In this section we prove Theorem~\ref{second}. Note that there are no ambiguous configurations on at most 3 points, so Theorem~\ref{second} clearly holds for $n \leq 3$.
Thus we have a deal with $n \geq 4$ to prove the theorem. First we deal with the planar case.

Let us denote by $\mathbb P_2^u\subset \mathbb P_2$ the subset of configurations having a unique Steiner tree. Observe that, due to Theorem~\ref{uniquepathconnectedness}, Theorem~\ref{second} will follow from the following
\begin{theo}
    Let $T_1, T_2$ be two embedding classes and $P_1,P_2\in \mathbb P_2^u$ be the two configurations of terminal points of the corresponding canonical realizations $\St_{EC_1},\St_{EC_2}$. Then there is a path in $\mathbb P_2^u$ connecting $P_1$ and $P_2$.
    \label{theo:path_between_topologies}
\end{theo}

Define the special (non-labelled) \textit{all-left linear tree} $ALT$ to be the path on $n$ vertices starting at $y_0$ and turning left at every branching point of $\Sigma_\infty$, i.e. 
$ALT$ is the subgraph of $\Sigma_\infty$ with the vertices $y_0,y_1, \dots, y_{2^k}, \dots, y_{2^{n-2}}$. To establish Theorem~\ref{theo:path_between_topologies} we will show that any $\St_{EC}$ corresponding to an embedding class with $n$ terminal vertices can be continuously deformed to $ALT$ inside the space of unique Steiner trees with some deformation preserving the labeling of terminal vertices. We construct such a deformation in several steps described below. In each step we continuously deform the set $P$ of the terminal points of $\St_{EC}$ to the set $P'$ of the terminals of $\St_{EC'}$ inside $\mathbb P_2^u$ by moving several points from $P$ one by one.

We need a preliminary lemma.

\begin{lemma}
\label{turning}
Let $EC$ be an embedding class with generic topology $R$. Suppose that $\St_{EC}$ contains a leaf $B = y_{2k}$  adjacent to the terminal $A = y_{k}$ of degree 2. 
Then one can continuously move $B$ to $y_{2k+1}$ along a path $\gamma$ in such a way that the whole configuration will remain in $\mathbb P_2^u$ at any point of $\gamma$ and has the embedding class $EC$.
\end{lemma}

\begin{proof}
 We construct a desired part of $\gamma$ explicitly (see Fig.~\ref{pic:lemma}).
First move $B$ into $\mu$-neighborhood of $A$ inside the segment $[y_{2k}y_{k}]$, where a small enough $\mu$ will be defined in the next paragraph; by Corollary~\ref{coruniquesubtree} a Steiner tree is unique and has the embedding class $EC$ at any configuration from this part of $\gamma$. 

\begin{figure}[h]
    \centering
    \begin{tikzpicture}[scale=0.8]

    \begin{scope}[shift={(-7,0)}]
        \draw[ultra thick, blue] (0,0) --++ (0,1) coordinate(A) --++ (150:2) coordinate (B);
        \draw[ultra thick, blue, dashed] (0,0) -- (0,-1);
        \draw ([yshift=-.2cm]A) node[xshift=-0.3cm, yshift=-0.2cm] {$\cfrac{2\pi}3$};
        \fill (A) circle (2pt) node[right]{$A = y_k$};
        \fill (B) circle (2pt) node[above]{$B = y_{2k}$};
    \end{scope}
    
    \begin{scope}[shift={(-3.5,0)}]
        \draw[ultra thick, blue] (0,0) --++ (0,1) coordinate(A) --++ (150:0.2) coordinate (B);
        \draw[ultra thick, blue, dashed] (0,0) -- (0,-1);
        \draw ([yshift=-.2cm]A) node[xshift=-0.3cm, yshift=-0.2cm] {$\cfrac{2\pi}3$};
        \fill (A) circle (2pt) node[right]{$A$};
        \fill (B) circle (2pt) node[above]{$B$};
        \draw[->] ($(A)!6!(B)$) arc (150:90:1);
    \end{scope}
    
    \begin{scope}[shift={(0,0)}]
        \draw[ultra thick, blue] (0,0) --++ (0,1) coordinate(A) --++ (90:0.2) coordinate (B);
        \draw[ultra thick, blue, dashed] (0,0) -- (0,-1);
        \draw ([yshift=-.2cm]A);
        \fill (A) circle (2pt) node[left]{$A$};
        \fill (B) circle (2pt) node[above]{$B$};
        \draw[->] ($(A)!6!(B)$) arc (90:30:1);
    \end{scope}
    
    \begin{scope}[shift={(3.5,0)}]
        \draw[ultra thick, blue] (0,0) --++ (0,1) coordinate(A) --++ (30:0.2) coordinate (B);
        \draw[ultra thick, blue, dashed] (0,0) -- (0,-1);
        \draw ([yshift=-.2cm]A) node[xshift=0.3cm, yshift=-0.2cm] {$\cfrac{2\pi}3$};
        \fill (A) circle (2pt) node[above left]{$A$};
        \fill (B) circle (2pt) node[above]{$B$};
    \end{scope}

    \begin{scope}[shift={(7,0)}]
        \draw[ultra thick, blue] (0,0) --++ (0,1) coordinate(A) --++ (30:2) coordinate (B);
        \draw[ultra thick, blue, dashed] (0,0) -- (0,-1);
        \draw ([yshift=-.2cm]A) node[xshift=0.3cm, yshift=-0.2cm] {$\cfrac{2\pi}3$};
        \fill (A) circle (2pt) node[above left]{$A$};
        \fill (B) circle (2pt) node[above]{$B = y_{2k+1}$};
    \end{scope}

\end{tikzpicture}
    \caption{The construction of $\gamma$ in Lemma~\ref{turning}.}
    \label{pic:lemma}
\end{figure}

Let $\bar R$ be the topology of $\St_{EC} \setminus [AB]$; obviously $\bar R$ is also generic.
Observation~\ref{obs:topologies}~(i) states that $\bar R$ lies in the exactly one set $D(\bar O)$, where $\bar O$ is a full topology.
By Corollary~\ref{coruniquesubtree} $\St_{EC} \setminus [AB]$ is a unique Steiner tree with $n-1$ terminals, so by 
Observation~\ref{obs:topologies}~(ii) there is $\eta > 0$ such that 
any other full topology $R'$ the length of the realization from $D(R')$ exceeds $\HH(\St \setminus [AB])$ by at least $\eta$.
Put $\mu = \eta/2$. 

Now rotate $B$ around $A$: let $B(\alpha)$, $\alpha \in [0,2\pi/3]$ be a such point that $|BA| = |B(\alpha)A|$ and the clockwise-oriented angle 
$y_{[k/2]}AB(\alpha)$ is equal to $2\pi/3 + \alpha$. In particular $B(0) \in [Ay_{2k})$, $B(2\pi/3) \in [Ay_{2k+1})$.

Let $\St(\alpha)$ be a Steiner tree  for the terminals of $\St \setminus [AB]$ and $B(\alpha)$.
Then
\[
\HH(\St(\alpha)) \leq  \HH(\St \setminus [AB] \cup [AB(\alpha)]) = \HH(\St \setminus [AB]) + \mu < \HH(\St \setminus [AB]) + \eta.
\]
Let $O$ be the full topology such that $R \in D(O)$.
Then the topology of $\St(\alpha)$ belongs to $D(O)$.
By Proposition~\ref{GPuniq} $\St(\alpha)$ is uniquely defined. The set $\St \setminus [AB] \cup [AB(\alpha)]$ is locally minimal, and has the topology from $D(R)$, so it coincides with $\St(\alpha)$. Thus not only the topology but the embedding class is preserved during this part of the path.

Finally, move $B$ from $B(2\pi/3)$ to $y_{2k}$ inside the segment $[y_{2k}y_{k}]$.
\end{proof}

%Let $T_1$ and $T_2$ be arbitrary topologies on $n$ terminals. 
%Let us connect their canonical realizations $\St_{T_1}$ and $\St_{T_2}$ by a path $\gamma \subset \mathbb{P}$ in a way that a Steiner tree is unique for every point of $\gamma$. 

%In the remainder of this section assume $i = 1,2$.

Let us now fix an embedding class $EC$ and construct the desired deformation of $\St_{EC}$ to $ALT$ inside the space of unique Steiner trees.

\paragraph{Step 1. Transform $\St_{EC}$ into a full Steiner tree $\St_{EC'}$ inside $\mathbb P_2^u$.} To make such a transformation we need to move all terminal vertices of $\St_{EC}$ of degree 2 or 3 to make them leaves.

Suppose first that $\St_{EC}$ contains a terminal $A = y_j$ of degree 2.
By the construction of $\St_{EC}$, the vertex $A$ is adjacent to vertices $B = y_{2j}$ and $C = y_{\lfloor j/2\rfloor}$, which may be terminals or Steiner points.
Move $A$ towards $y_{2j+1}$ along the edge $y_jy_{2j+1}$ of $\Sigma_\infty$ until it hits $y_{2j+1}$ (see. Fig.~\ref{nodegree2}). 
\begin{figure}[h]
    \begin{center}
    \begin{tikzpicture}[scale=1.2]

    \begin{scope}[shift={(-3.5,0)}]
        \draw[ultra thick, blue] (0,0) --++ (0,1) coordinate(A) --++ (150:2) coordinate (B);
        \draw[ultra thick, blue, dashed] (B) --++ (210:0.7);
        \draw[ultra thick, blue, dashed] (B) --++ (90:0.7);
        \draw[ultra thick, blue, dashed] (0,0) -- (0,-1) coordinate(C);
        \draw[ultra thick, blue, dashed] (0,0) -- (0,-1);
        \draw ([yshift=-.2cm]A) node[left] {$\cfrac{2\pi}3$};
        \fill (A) circle (2pt) node[right]{$A = y_{j}$};
        \fill (B) circle (2pt) node[right]{$B = y_{2j}$};
        \fill (C) circle (0pt) node[right]{$C = y_{[j/2]}$};
    \end{scope}
    
    \begin{scope}[shift={(0,0)}]
        \draw[ultra thick, blue] (0,0) --++ (0,1) coordinate(Y) --++ (150:2) coordinate (B);
        \draw[ultra thick, blue] (Y) --++ (30:2) coordinate(A);
        \draw[ultra thick, blue, dashed] (B) --++ (210:0.7);
        \draw[ultra thick, blue, dashed] (B) --++ (90:0.7);
        \draw[ultra thick, blue, dashed] (0,0) -- (0,-1) coordinate(C);
        \draw[ultra thick, blue, dashed] (0,0) -- (0,-1);
        \fill (A) circle (2pt) node[above]{$A = y_{2j+1}$};
        \fill (B) circle (2pt) node[right]{$B = y_{2j}$};
        \fill (C) circle (0pt) node[right]{$C = y_{[j/2]}$};
    \end{scope}

\end{tikzpicture}
    \caption{Elimination of points with degree 2 in $\St_{EC}$}
    \label{nodegree2}
    \end{center}
\end{figure}
Corollary~\ref{coruniquesubtree} ensures that thus obtained deformation of the set of terminal points lies inside $\mathbb P_2^u$. Applying this deformation to each terminal vertex of degree 2 one by one we eventually get rid of those.

Assume now that $\St_{EC}$ has a terminal point of degree $3$. Since $\St_{EC}$ has no terminal of degree two and the number of Steiner points is at most the number of leafs minus two, one may move terminals of degree three one by one in a neighborhoods of different leaves by a path in $\mathbb{P}$. From now a topology of a tree is generic.

Now consider any point $A$ in an $\eps$-neighborhood of a leaf $B$ for some small $\eps$. Then continue moving $A$ while moving $B$ simultaneously in the same direction until $A$ riches $y_k$ and $B$ riches $B(\pi/3)$ (see Fig.~\ref{nodegree3}). Now stop moving $A$, but rotate $B$ around $A$ until it hits the ray $Ay_{2k+1}$, then extend $B$ to $y_{2k+1}$ and $A$ to $y_{2k}$. Now all our terminal points again belongs to the set $\{y_0,y_1,\dots\}$ and the unique Steiner tree is given by the canonical realization $\St_{EC'}$ for some new embedding class $EC'$. The fact that the set of terminal points was staying inside $\mathbb P_2^u$ while we were moving them follows from Corollary~\ref{coruniquesubtree} and the proof of Lemma~\ref{turning}.

\begin{figure}[h]
    \begin{center}
    \begin{tikzpicture}[scale=1.2]

    \begin{scope}[shift={(0,0)}]
        \draw[ultra thick, blue] (0,0) --++ (0,0.5) coordinate(A) --++ (0,0.5) coordinate (B);
        \draw[ultra thick, blue, dashed] (0,0) -- (0,-0.7);
        \fill (A) circle (2pt) node[above left]{$A$};
        \fill (B) circle (2pt) node[above]{$B = y_k$};
    \end{scope}t={(1,0)}]
    
    \begin{scope}[shift={(2,0)}]
         \draw[ultra thick, blue] (0,0) --++ (0,1) coordinate(A) --++ (0,0.5) coordinate (B);
        \draw[ultra thick, blue, dashed] (0,0) -- (0,-0.7);
        \fill (A) circle (2pt) node[above left]{$A = y_k$};
        \fill (B) circle (2pt) node[above]{$B = B(\pi/3)$};
        \draw[->] ($(A)!0.6!(B)$) arc (90:30:0.5);
    \end{scope}

    \begin{scope}[shift={(4.5,0)}]
         \draw[ultra thick, blue] (0,0) --++ (0,1) coordinate(A) --++ (30:0.5) coordinate (B);
        \draw[ultra thick, blue, dashed] (0,0) -- (0,-0.7);
        \fill (A) circle (2pt) node[above left]{$A = y_k$};
        \fill (B) circle (2pt) node[above]{$B = B(\pi/3)$};
    \end{scope}
    
     \begin{scope}[shift={(8,0)}]
         \draw[ultra thick, blue] (0,0) --++ (0,1) coordinate(Y) --++ (30:1) coordinate (B);
          \draw[ultra thick, blue] (Y) --++ (150:1) coordinate (A);
        \draw[ultra thick, blue, dashed] (0,0) -- (0,-0.7);
        \fill (A) circle (2pt) node[above left]{$A = y_{2k}$};
        \fill (B) circle (2pt) node[above]{$B = y_{2k+1}$};
    \end{scope}

    %\foreach\x in {-2.5,0,2.5}{
    %    \draw (\x,1) node{$\longrightarrow$};
    %}
\end{tikzpicture}
    \caption{Elimination of points with degree 3 in $\St_{EC}$}
    \label{nodegree3}
    \end{center}
\end{figure}

Note that $\St_{EC'}$ still has no terminal points of degree 2 and has one less terminal point of degree 3 than $\St_{EC}$. Hence we can do this procedure until we obtain a canonically realized full tree. 

%A continuous movement of $A$ along the relative interior of $S$ corresponds to a path in $\mathbb{P}$ since $S$ has no terminal points in the relative interior.
%Move $A$ into a small neighborhood of an arbitrary terminal vertex $B = y_k$ of $S$; by the definition $y_k$ has degree 1.
%Then move $B$ in $B(\pi/3)$ (in the notation of Lemma~\ref{turning}) and $A$ to $y_k$. 
%Finally, we rotate $B$ to $B(0) = y_{2k}$ and then move $A$ to $y_{2k+1}$. By Corollary~\ref{coruniquesubtree} and Lemma~\ref{turning} Steiner tree is unique during the path.

%The resulting tree $\St_T$ is full.

\paragraph{Step 2. Permute the labels of terminal points of $\St_{EC}$ if necessary.}
Now we can assume that $\St_{EC}$ is a full tree. By Theorem~\ref{uniquepathconnectedness} one may put label 1 into the root by a path of path in $\mathbb{P}$ and do not touch the root of the tree later on. 
Let $A$ and $B$ be two terminal points which we want to swap.

Now we can assume that $\St_{EC}$ is a full tree, hence it has exactly $n-2$ Steiner points. In particular, we can choose two Steiner points of $\St_{EC}$ that are connected with two terminal points of $\St_{EC}$; let $y_k$ be the one of them which is not adjacent to the root of $\Sigma_\infty$. Denote the terminals adjacent to $y_k$ by $B = y_{2k}$ and $A = y_{2k+1}$.

We may swap $A$ with any label. First swap $A$ and $B$ as shown at Fig.~\ref{swapneighbors}: 
move $B$ into $y_k$ and $A$ in a small neighborhood of $y_k$, then turn and finally make a reverse procedure.
By Lemma~\ref{turning} Steiner tree is unique during the middle part of this procedure; by Corollary~\ref{coruniquesubtree} Steiner tree is unique during other parts.

\begin{figure}[h]
    \begin{center}
    \begin{tikzpicture}[scale=0.82]

    \begin{scope}[shift={(-8,0)}]
        \draw[ultra thick, blue] (0,0) --++ (0,1) coordinate (b) --++ (150:2) coordinate (B);
        \draw[ultra thick, blue] (0,0) --++ (0,1) coordinate (b) --++ (30:2) coordinate (A);
        \draw[ultra thick, blue, dashed] (0,0) -- (0,-1);
        \draw ([yshift=-.2cm]b) node[xshift=-0.3cm, yshift=-0.2cm] {$\cfrac{2\pi}3$};
        \fill (A) circle (2pt) node[above]{$A = y_{2k+1}$};
        \fill (B) circle (2pt) node[above]{$B = y_{2k}$};
    \end{scope}

    \begin{scope}[shift={(-3.5,0)}]
        \draw[ultra thick, blue] (0,0) --++ (0,1) coordinate(A) --++ (150:2) coordinate (B);
        \draw[ultra thick, blue, dashed] (0,0) -- (0,-1);
        \draw ([yshift=-.2cm]A) node[xshift=-0.3cm, yshift=-0.2cm] {$\cfrac{2\pi}3$};
        \fill (A) circle (2pt) node[right]{$A = y_{k}$};
        \fill (B) circle (2pt) node[below]{$B$};
    \end{scope}
    
    \begin{scope}[shift={(-1,0)}]
        \draw[ultra thick, blue] (0,0) --++ (0,1) coordinate(A) --++ (150:1) coordinate (B);
        \draw[ultra thick, blue, dashed] (0,0) -- (0,-1);
        \draw ([yshift=-.2cm]A) node[xshift=-0.3cm, yshift=-0.2cm] {$\cfrac{2\pi}3$};
        \fill (A) circle (2pt) node[right]{$A$};
        \fill (B) circle (2pt) node[above]{$B$};
        \draw[->] ($(A)!.6!(B)$) arc (150:30:.6);
    \end{scope}
    
    \begin{scope}[shift={(1,0)}]
        \draw[ultra thick, blue] (0,0) --++ (0,1) coordinate(A) --++ (30:1) coordinate (B);
        \draw[ultra thick, blue, dashed] (0,0) -- (0,-1);
        \draw ([yshift=-.2cm]A) node[xshift=0.3cm, yshift=-0.2cm] {$\cfrac{2\pi}3$};
        \fill (A) circle (2pt) node[above left]{$A$};
        \fill (B) circle (2pt) node[above]{$B$};
    \end{scope}t={(1,0)}]
    
    \begin{scope}[shift={(3.5,0)}]
        \draw[ultra thick, blue] (0,0) --++ (0,1) coordinate(A) --++ (30:2) coordinate (B);
        \draw[ultra thick, blue, dashed] (0,0) -- (0,-1);
        \draw ([yshift=-.2cm]A) node[xshift=0.3cm, yshift=-0.2cm] {$\cfrac{2\pi}3$};
        \fill (A) circle (2pt) node[above left]{$A$};
        \fill (B) circle (2pt) node[above]{$B$};
    \end{scope}
    
    \begin{scope}[shift={(8,0)}]
        \draw[ultra thick, blue] (0,0) --++ (0,1) coordinate (b) --++ (150:2) coordinate (A);
        \draw[ultra thick, blue] (0,0) --++ (0,1) coordinate (b) --++ (30:2) coordinate (B);
        \draw[ultra thick, blue, dashed] (0,0) -- (0,-1);
        \draw ([yshift=-.2cm]b) node[xshift=-0.3cm, yshift=-0.2cm] {$\cfrac{2\pi}3$};
        \fill (A) circle (2pt) node[above]{$A$};
        \fill (B) circle (2pt) node[above]{$B$};
    \end{scope}

  %  \foreach\x in {-5.8,-2.5,0,2.5,5.8}{
  %      \draw (\x,1) node{$\longrightarrow$};
  %   }
\end{tikzpicture}
    \caption{Swapping the labels of terminals connecting with a common branching point}
    \label{swapneighbors}
    \end{center}
\end{figure}

Then swap $A$ with any terminal $C \neq B$ of $\St_i'$ (see Fig.~\ref{swapfar}).
Start with the previous procedure and stop it at the point $B = B(\pi/3)$ (in the notation of Lemma~\ref{turning}). Then $A$ moves inside the tree into a neighborhood of $C = y_l$ and $B$ comes to $y_k$. 
We are going to apply Lemma~\ref{turning} to $A$ and $C$: move $C$ to $C(\pi/3)$ and $A$ to $y_l$. 
Then $C$ rotates to $y_{2l+1}$, after that $A$ moves to $y_{2l}$. 
Now the positions of $A$ and $C$ are symmetric so we may do the reverse procedure after swapping $A$ and $C$.

\begin{figure}[h]
    \begin{center}
    \begin{tikzpicture}[scale=1]

    \begin{scope}[shift={(-3,0)}]
        \draw[ultra thick, blue] (0,0) --++ (0,1) coordinate(Y) --++ (150:1) coordinate (B);
        \draw[ultra thick, blue] (Y) --++ (30:1) coordinate (A);
        \draw[ultra thick, blue, dashed] (0,0) --++ (0,-1)coordinate (X) --++ (-30:1) coordinate (Z) --++ (30:0.5) coordinate (T);
        \draw[ultra thick, blue, dashed] (X) --++ (-150:0.5);
        \draw[ultra thick, blue] (T)--++ (30:0.5) coordinate (C);
        
        \fill (A) circle (2pt) node[above]{$A=y_{2k+1}$};
        \fill (B) circle (2pt) node[above]{$B=y_{2k}$};
        \fill (C) circle (2pt) node[above]{$C=y_{l}$};
    \end{scope}
    
    \begin{scope}[shift={(0,0)}]
        \draw[ultra thick, blue] (0,0) --++ (0,1) coordinate(A) --++ (150:0.5) coordinate (B);
        \draw[ultra thick, blue, dashed] (0,0) -- (0,-1);
        \draw[->] ($(A)!.6!(B)$) arc (150:110:.6);
        \draw[ultra thick, blue, dashed] (0,0) --++ (0,-1)coordinate (X) --++ (-30:1) coordinate (Z) --++ (30:0.5) coordinate (T);
        \draw[ultra thick, blue, dashed] (X) --++ (-150:0.5);
        \draw[ultra thick, blue] (T)--++ (30:0.5) coordinate (C);
        
        \fill (A) circle (2pt) node[right]{$A = y_k$};
        \fill (B) circle (2pt) node[above]{$B$};
        \fill (C) circle (2pt) node[right]{$C$};
    \end{scope}
    
    \begin{scope}[shift={(3,0)}]
        \draw[ultra thick, blue] (0,1)coordinate(B) --++ (0,-0.5) coordinate(A);
        \draw[ultra thick, blue, dashed] (A) -- (0,-1);
        
         \draw[ultra thick, blue, dashed] (0,0) --++ (0,-1)coordinate (X) --++ (-30:1) coordinate (Z) --++ (30:0.5) coordinate (T);
        \draw[ultra thick, blue, dashed] (X) --++ (-150:0.5);
        \draw[ultra thick, blue] (T)--++ (30:0.5) coordinate (C);
        
        \fill (A) circle (2pt) node[right]{$A$};
        \fill (B) circle (2pt) node[above]{$B = y_k$};
        \fill (C) circle (2pt) node[right]{$C$};
    \end{scope}

    \begin{scope}[shift={(6,0)}]
        \draw[ultra thick, blue] (0,1)coordinate(B) --++ (0,-0.5) coordinate(A);
        \draw[ultra thick, blue, dashed] (A) -- (0,-1);
        
         \draw[ultra thick, blue, dashed] (0,0) --++ (0,-1)coordinate (X) --++ (-30:1) coordinate (Z) --++ (30:0.5) coordinate (T);
        \draw[ultra thick, blue, dashed] (X) --++ (-150:0.5);
        \draw[ultra thick, blue] (T)--++ (30:0.5) coordinate (C);
        
        \fill (T) circle (2pt) node[above]{$A$};
        \fill (B) circle (2pt) node[above]{$B$};
        \fill (C) circle (2pt) node[right]{$C$};
    \end{scope}

    \begin{scope}[shift={(9,0)}]
        \draw[ultra thick, blue] (0,1)coordinate(B) --++ (0,-0.5) coordinate(A);
        \draw[ultra thick, blue, dashed] (A) -- (0,-1);
        
         \draw[ultra thick, blue, dashed] (0,0) --++ (0,-1)coordinate (X) --++ (-30:1) coordinate (Z) --++ (30:0.5) coordinate (T);
        \draw[ultra thick, blue, dashed] (X) --++ (-150:0.5);
        \draw[ultra thick, blue] (T)--++ (30:0.5) coordinate (C) --++ (-30:0.5) coordinate(CC);
        \draw[ultra thick, blue] (C)--++ (90:0.5) coordinate (AA);
        
        \fill (AA) circle (2pt) node[above]{$A = y_{2l}$};
        \fill (B) circle (2pt) node[above]{$B = y_k$};
        \fill (CC) circle (2pt) node[below]{$C = y_{2l+1}$};
    \end{scope}
    
\end{tikzpicture}
    \caption{Swapping the labels of arbitrarily terminals}
    \label{swapfar}
    \end{center}
\end{figure}

Finally to swap labels of arbitrary terminals $C$ and $D$ we swap $A = y_{2k+1}$ with $C$, $C = y_{2k+1}$ with $D$ and $D  = y_{2k+1}$ with $A$.
Since the set of all transpositions spans the symmetric group we may construct a path in $\mathbb{P}$ connecting 
$\St_{EC}$ and the same tree with an arbitrary permutation of its labels.
Till the end of the section all trees are not labelled.

\paragraph{Step 3. Connect $\St_i'$ with the all-left linear tree $ALT$ by a path in $\mathbb{P}$.} 
While there is a terminal point $A = y_j$ of $\St_i'$ not belonging to $ALT$, consider such a vertex with the largest $j$.
It implies that the degree of $A$ is 1. Our aim is to move $A$ inside $\Sigma_\infty$ to the first vertex $y_{w}$ of $ALT$ which does not belong to $\St_i'$.

Consider the case when $A$ is adjacent to a branching point $y_l$ then $j = 2l+1$ and $B = y_{2l}$ is also a terminal of $\St_i'$ because of the maximality of $j$. Move $A$ into $y_l$ and rotate $B$ into $B(\pi/3)$ (in the notation of Lemma~\ref{turning}).
Then $A$ moves into the tree and $B$ moves into $y_l$. 

Now $A$ is either inside the tree or $A$ is a terminal connected with a vertex of degree 2.
Move $A$ into $y_w$, the only problem is that $A$ cannot coincide with the terminal of degree 2.
Movement through a terminal of degree 2 is depicted in Fig.~\ref{arounddegree2}.

\begin{figure}[h]
    \begin{center}
    \begin{tikzpicture}[scale=1.2]

    \begin{scope}[shift={(-3.5,0)}]
        \draw[ultra thick, blue] (0,0.5) coordinate(A) --++ (0,1) coordinate(Y) --++ (150:1.5) coordinate (B);
        \draw[ultra thick, blue, dashed] (0,0.5) -- (0,-0.5);
        \draw ([yshift=-.2cm]Y) node[xshift=-0.3cm, yshift=-0.2cm] {$\cfrac{2\pi}3$};
        \draw[ultra thick, blue, dashed] (B) -- ($(B)!-0.8!(Y)$);

        \fill (A) circle (2pt) node[right]{$A$};
        \fill (Y) circle (2pt) node[xshift=0.1cm, yshift=0.5cm]{$y_k = B$};

    \end{scope}

    \begin{scope}[shift={(0,0)}]
        \draw[ultra thick, blue] (0,0.5) coordinate(A) --++ (0,1) coordinate(Y) --++ (150:1.5) coordinate (B);
        \draw[ultra thick, blue, dashed] (0,0.5) -- (0,-0.5);
        \draw[ultra thick, blue] (Y) --++ (30:0.5) coordinate(BB);

        \draw ([yshift=-.2cm]Y) node[xshift=-0.3cm, yshift=-0.2cm] {$\cfrac{2\pi}3$};
        \draw[ultra thick, blue, dashed] (B) -- ($(B)!-0.8!(Y)$);

        \fill (A) circle (2pt) node[right]{$A$};
        \fill (Y) circle (2pt) node[xshift=0.0cm, yshift=0.4cm]{$y_k$};
        \fill (BB) circle (2pt) node[right]{$B$};

    \end{scope}
    
    \begin{scope}[shift={(3.5,0)}]
        \draw[ultra thick, blue] (0,0.5) coordinate(A) --++ (0,1) coordinate(Y) --++ (150:1.5) coordinate (B);
        \draw[ultra thick, blue, dashed] (0,0.5) -- (0,-0.5);
        \draw[ultra thick, blue] (Y) --++ (30:0.5) coordinate(BB);

        \draw ([yshift=-.2cm]Y) node[xshift=-0.3cm, yshift=-0.2cm] {$\cfrac{2\pi}3$};
        \draw[ultra thick, blue, dashed] (B) -- ($(B)!-0.8!(Y)$);

        \fill (B) circle (2pt) node[xshift=0.1cm, yshift=0.3cm]{$A$};
        \fill (Y) circle (2pt) node[xshift=0.0cm, yshift=0.4cm]{$y_k$};
        \fill (BB) circle (2pt) node[right]{$B$};
    \end{scope}

    \begin{scope}[shift={(7,0)}]
        \draw[ultra thick, blue] (0,0.5) coordinate(A) --++ (0,1) coordinate(Y) --++ (150:1.5) coordinate (B);
        \draw[ultra thick, blue, dashed] (0,0.5) -- (0,-0.5);
        \draw ([yshift=-.2cm]Y) node[xshift=-0.3cm, yshift=-0.2cm] {$\cfrac{2\pi}3$};
        \draw[ultra thick, blue, dashed] (B) -- ($(B)!-0.8!(Y)$);

        \fill (B) circle (2pt) node[xshift=0.1cm, yshift=0.3cm]{$A$};
        \fill (Y) circle (2pt) node[xshift=0.1cm, yshift=0.5cm]{$y_k = B$};

    \end{scope}

\end{tikzpicture}
    \caption{Movement through a terminal of degree 2}
    \label{arounddegree2}
    \end{center}
\end{figure}

Finally all the vertices of $\St_i'$ belong to $ALT$, so we are done.
Since we connect $\St_{EC_1}$ with $ALT$ and $\St_{EC_2}$ with $ALT$, the desired $\gamma$ is constructed.

\begin{proof}[Proof of Theorem~\ref{second}]
Let $P_1, P_2 \in \mathbb{P}_d$ be configurations with unique Steiner trees $\St(P_1)$ and $\St(P_2)$ having embedding classes $T_1$ and $T_2$, respectively.
By Theorem~\ref{uniquepathconnectedness} there is a path $\gamma_i$ between $\St(P_i)$ and $\St_{T_i}$ in $\mathbb{P}_d$ such that a Steiner tree is unique during $\gamma_i$. 
We have constructed the path $\gamma$ between $\St_{T_1}$ and $\St_{T_2}$ in $\mathbb{P}_2 \subset \mathbb{P}_d$; a Steiner tree is also unique during $\gamma$.
The gluing of $\gamma_1$, $\gamma$ and $\gamma_2^{-1}$ finishes the proof.
\end{proof}

\section{Proof of Theorem~\ref{main}}
\label{sec:proof of main theorem}

This section is devoted to the proof of Theorem~\ref{main}. The prove is using the theory of subanalytic sets, and for the sake of completeness we begin our exposition with a brief reminder of some definitions and facts from this theory.

\subsection{Subanalytic subsets of a real analytic manifold}
\label{subsec:Subanalytic subsets of a real analytic manifold}
All the facts expounded in this section are well-known. During our exposition we mostly follow Sections~2 and~3 from the paper~\cite{bierstone1988semianalytic}.

Let $M$ be a real analytic manifold and $\Oo_M$ denote the sheaf of real analytic functions on $M$, that is, for any open $U\subset M$ the set $\Oo_M(U)$ is the space of real analytic functions defined on $U$. We introduce the following definitions:
\begin{enumerate}
    \item A subset $A\subset M$ is called an \emph{analytic submanifold} if for any $p\in A$ there exists a neighborhood $U\subset M$ of $p$ in $A$ such that either $A\cap U = U$, or there exist a finite collection $f_1,f_2,\dots,f_k\in \Oo_M(U)$ such that $A\cap U$ is the set of common zeros of $f_1,\dots, f_k$ and for any $x\in A\cap U$ the gradients $\nabla f_1(x),\dots,\nabla f_k(x)$ are linearly independent.
    \item A subset $A\subset M$ is called~\emph{analytic} if for any $p\in M$ there exists a neighborhood $U$ of $p$ such that either $A\cap U = U$, or there exists a finite set of functions $f_1,\dots, f_k\in \Oo_M(U)$ such that $A\cap U$ is the set of common zeros of $f_1,\dots,f_k$. Note that we require this property for all $p\in M$, not only for $p\in A$.
    \item A subset $A\subset M$ is called \emph{semianalytic} if for any point $p\in M$ there exists a neighborhood $U\subset M$ and a finite number of subsets $A_{i,j}\subset U$ such that $A\cap U = \cup_i\cap_j A_{i,j}$ and each $A_{i,j}$ is of the form $\{f>0\}$ or $\{f=0\}$ for some $f\in \Oo_M(U)$. A semianalytic subset $A$ is called \emph{smooth} if it is an analytic submanifold.
\end{enumerate}

The following lemma follows from~\cite[Proposition~2.10]{bierstone1988semianalytic}:
\begin{lemma}
    \label{lemma:semianalytic_local_structure}
    Let $M$ be a real analytic manifold, $A\subset M$ be a semianalytic subset and $p\in M$ be an arbitrary point. Then there exists a neighborhood $U\subset M$ of $p$ and a finite collection of disjoint subsets $A_1,A_2,\dots,A_k\subset U$ such that
    \begin{enumerate}
        \item each of $A_1,\dots, A_k$ is a semianalytic subset of $U$ and an analytic submanifold of $M$, and
        \item $A\cap U$ is a disjoint union of $A_1,\dots, A_k$.
    \end{enumerate}
\end{lemma}

Semianalytic sets admit many properties similar to those of semialgebraic sets (i.e. those given by polynomial inequalities), but the theories are not identical. An important difference is that images of semianalytic sets under proper maps are not necessary semianalytic (see~\cite[Example~2.14]{bierstone1988semianalytic}), while for semialgebraic sets this is always true. This motivates the following definition:

\begin{definition}
    A subset $X\subset M$ is called subanalytic if for any $p\in M$ there is a neighborhood $U$ of $p$, an analytic manifold $N$ and a relatively compact semianalytic subset $A\subset M\times N$ such that $X\cap U = \pi(A)$, where $\pi$ is the projection on $M$.
\end{definition}

The following lemma follows immediately from this definition:
\begin{lemma}
    \label{lemma:images_of_semianalytic_are_subanalytic}
    Let $N,M$ be two analytic manifolds and $f: N\to M$ be an analytic map. Assume that $A\subset N$ is semianalytic and $f$ restricted to the closure of $A$ is proper. Then $f(A)$ is a subanalytic subset of $M$.
\end{lemma}
\begin{proof}
    Let $\Gamma_f\subset M\times N$ be the graph of the mapping $f$. The graph is an analytic subset of $M\times N$ since $f$ is analytic. Let $p\in M$ and $U\subset M$ be its relatively compact semianalytic neighborhood. Define $B = (U\times A)\cap \Gamma_f\subset M\times N$, then $B$ is semianalytic and relatively compact since since $f$ restricted to the closure of $A$ is proper. We have $f(A)\cap U = \pi(B)$, where $\pi$ is the projection on $M$. Since $p$ were arbitrary, we conclude that $f(A)$ is subanalytic.
\end{proof}

Subanalytic sets, although not being semianalytic in general, still have a lot of nice properties. Direct products, finite intersections and unions, closures, complements and, thus, interiors of subanalytic sets are still subanalytic (see~\cite[Chapter~3]{bierstone1988semianalytic}). The following lemma describes a local structure of subanalytic sets, see~\cite[Lemma~3.4]{bierstone1988semianalytic}):

\begin{lemma}
    \label{lemma:local_structure_of_subanalytic_sets}
    Let $N,M$ be analytic manifolds and $A\subset N\times M$ be a relatively compact semianalytic subset. Then there exists a finite collection of smooth connected semianalytic subsets $A_1,\dots,A_k\subset N\times M$ such that
    \begin{enumerate}
        \item $A = \sqcup_{j = 1}^k A_j$,
        \item for any $j$ the rank of $d\pi$ on $T_xA_j$ does not depend on $x\in A_j$.
    \end{enumerate}
\end{lemma}

From this lemma we get an immediate corollary:

\begin{cor}
    \label{cor:local_structure_of_subanalytic_sets}
    Let $M$ be an analytic manifold and $X\subset M$ be a subanalytic subset. Theb there exists a countable collection of connected analytic submanifolds $X_1,X_2,X_3,\dots$ of $M$ such that $X = X_1\cup X_2\cup X_3\cup\ldots$.
\end{cor}
\begin{proof}
    Since the topology of $M$ has a countable base, it is enough to prove the statement of the corollary locally. Passing to a neighborhood of some point if necessary we can assume that there is a relatively compact semianalytic subset $A\subset M\times N$ for some real analytic manifold $N$ such that $X = \pi(A)$. Let $A_1,\dots,A_k\subset M\times N$ be such in Lemma~\ref{lemma:local_structure_of_subanalytic_sets}. For each $j = 1,\dots, k$ there is a countable collection of open subsets $U_{j1},U_{j2},\dots$ of $M\times N$ covering $A_j$ and such that $\pi(U_{ji}\cap A_j)$ is a connected analytic submanifold of $M$. Then we have
    \[
        X = \pi(A) = \bigcup_{j = 1}^k \bigcup_{i\geq 1} \pi(U_{ji}\cap A_j).
    \]
\end{proof}

For a technical reason we need to introduce the notion of a fiber product. Let $X,Y,U$ be some sets and $f:X\to U$, $g: Y\to U$ be some maps between these sets. Then
\begin{equation}
    \label{eq:def_of_fiber_product}
    X\times_{f=g} Y = \{(x,y)\in X\times Y\ \mid\ f(x) = g(y)\}.
\end{equation}
Note that we have a natural projection $X\times_{f=g} Y\to U$ which sends $(x,y)$ to $f(x)$.

\begin{lemma}
    \label{lemma:fiber_product}
    Assume that $M,N,U$ are real analytic manifolds and $f: M\to U$ and $g: N\to U$ are real analytic maps. Let$X\subset M$ and $Y\subset N$ be sub- or semianalytic subsets. Then $X\times_{f=g}Y$ is a sub- or semianalytic subset of $M\times_{f=g}N$ respectively.
\end{lemma}
\begin{proof}
    Follows immediately from definitions. Indeed, $X\times_{f=g}Y$ is the intersection of the sub- or semianalytic set $X\times Y$ and the subset $\{(x,y)\in M\times N\ \mid\ f(x) = g(y)\}$ inside $M\times N$, hence is sub- or semianalytic respectively.
\end{proof}

\subsection{Proof of Theorem~\ref{main}}
\label{subsec:proof_of_Theorem_on_dimension}

In this section we prove Theorem~\ref{main}. Recall that a number $n\geq 4$ of terminals is fixed and $\PP_2$ is equal to $(\mathbb R^2)^n$ with diagonals removed. Denote the subset of ambiguous configurations by $\Aa$. Let also $\Aa_{\mathrm{non-generic}}$ denote the set of all configurations admitting non-generic Steiner tree. Recall that $\dim \Aa_{\mathrm{non-generic}} = 2n-2$ by Observation~\ref{obs:topologies}.

We begin with the following
\begin{lemma}
    \label{lemma:dim_of_Amb}
    The Hausdorff dimension of $\Aa$ is at least $2n-1$.
\end{lemma}
\begin{proof}
    Let $T_1,\dots,T_N$ be all possible full topologies on $n$ points and $V(T_j)$ denotes the set of vertices of $T_j$. Given a map $f:V(T_j)\to \mathbb R^2$, let $L(f)$ the total length of the segments connecting $f(V(T_j))$ accordingly to topology $T_j$ (note that we do not claim any restrictions, in particular absence of cycles or local minimality), i.e.
    \[
        L(f) = \sum_{vw\text{ is an edge of }T_j} |f(v) - f(w)|.
    \]
    As it follows from Proposition~\ref{GPuniq}, for any $P\in \PP_2$ and $j = 1,\dots, N$ there exists precisely one map $f: V(T_j)\to \mathbb R^2$ which maps the terminals of $T_j$ to the points from $P$ keeping the enumeration and which minimize $L(f)$ among all such maps. Set $L_j(P) = L(f)$ in this case. Note that $L_j$ is a continuous function on $\PP_2$. Define
    \[
        B_j = \{ P\in \PP_2\ \mid\ L_j(P) < L_i(P)\ \forall\ i\neq j \}.
    \]
    It follows that $B_j$'s are open and disjoint sets. We also have $B_j\neq \varnothing$; indeed, by Corollary~\ref{coruniquesubtree} each $T_j$ is the topology of some Steiner tree which is unique. Note that $A = \PP_2 \smallsetminus \left(\cup_{j = 1}^N B_j\right)\subset \Aa\cup \Aa_{\mathrm{non-generic}}$ by Observation~\ref{obs:topologies}.
    The lemma now follows from $\dim \Aa_{\mathrm{non-generic}} = \dim\mathbb{R}^{2n}\smallsetminus \PP_2 = 2n-2$ and Lemma~\ref{lemma:dimension_of_complement}.
\end{proof}

\begin{lemma}
    \label{lemma:dimension_of_complement}
    Assume that $N \geq 2$, $m\geq 1$ and $B_1,\dots,B_N\subset \mathbb R^m$ are non-empty disjoint open sets. Put $A = \mathbb{R}^m\smallsetminus \left(\cup_{j = 1}^N B_j\right)$. Then $\dim A\geq m-1$, where $\dim$ is the Hausdorff dimension.
\end{lemma}
\begin{proof}
    Let $P_1,P_2\in \mathbb R^m$ be such that $P_1\in B_1$ and $P_2\in B_2$. Let $l$ be the line passing through these points and $V\subset \mathbb R^m$ be the subspace of codimension 1 orthogonal to $l$. Let $\pi: \mathbb R^m\to V$ be the orthogonal projection. Note that $v\in \pi(A)$ if and only if the line $\pi^{-1}(v)$ intersects $A$. In particular, $v_0 = \pi(l)\in \pi(A)$ and moreover there exists $\eps>0$ such that $v\in \pi(A)$ if $|v-v_0|\leq \eps$ since $B_1,B_2$ are open. It follows that $\pi(A)$ has a non-empty interior as a subset of $V$ and $\dim \pi(A) = m-1$. Since $\pi$ is 1-Lipshitz, it implies that $\dim A\geq m-1$.
\end{proof}

The converse estimate $\dim \Aa \leq 2n-1$ is more involved and requires some additional constructions. Let $T$ be some (not necessary full or generic) topology with $n$ terminals. Enumerate the Steiner points of $T$ arbitrary, let $k$ be the total amount of them. Given $(P,q) = (p_1,\dots, p_n,q_1,\dots, q_k)\in \PP_2 \times (\mathbb R^2)^k$, we can idendify the corresponding points on the plane with the vertices of $T$ following the enumeration and connect a pair of corresponding points by a straight segment for each edge of $T$. Thus, any such $(P,q)$ defines a map from $T$ to the plane. Let $\Rg(T)\subset \PP_2 \times (\mathbb R^2)^k$ be the following:
\[
    \Rg(T) = \{ (P,q)\in \PP_2 \times (\mathbb R^2)^k\ \mid\ (P,q)\text{ defines an embedding of $T$} \}.
\]
Obviously, $\Rg(T)$ is an open subset of $\PP_2 \times (\mathbb R^2)^k$. Let $L_T:  \PP_2 \times (\mathbb R^2)^k\to \mathbb R$ be the function that computes the length of the image of $T$; note that $L_T$ is real analytic on $\Rg(T)$ and continuous everywhere. Define
\[
    \Rc(T) = \{ (P,q)\in \Rg(T)\ \mid\ \nabla_q L_T(P,q) = 0 \},
\]
where $\nabla_q$ is the gradient with respect to the variable $q$. Since $L_T$ is real-analytic, $\Rc(T)$ is an analytic subset of $\Rg(T)$. We have the following

\begin{lemma}
    \label{lemma:properties_of_Rc}
    The following statements hold:
    \begin{enumerate}
        \item Let $(P,q)\in \Rc(T)$. Consider the function $L_T(P,\cdot)$ as a continuous function from $(R^2)^k$ to $\mathbb R$. Then $q$ is the unique point of the global minimum of $L_T(P,\cdot)$.
        \item Let $P_T:\Rg(T)\to \PP_2$ be the projection. Then $P_T$ restricted to $\Rc(T)$ is injective and the set $P_T(\Rc(T))$ is open in $\PP_2$.
    \end{enumerate}
\end{lemma}

\begin{proof}
    Note that for any fixed $P$ the function $L_T(P,\cdot)$ is coercive on $(\mathbb R^2)^k$. It follows that for any $P$ there is a point $q(P)\in (\mathbb R^2)^k$ where $L_T(P,\cdot)$ attains its global minimum. By Proposition~\ref{GPuniq} such a point is always unique and there are no other local minima of $L_T(P,\cdot)$. In particular, the mapping $P\mapsto q(P)$ is a well-defined mapping $\PP_2\to (\mathbb R^2)^k$. From the continuity of $L_T$ it follows that this mapping is continuous.

    To prove the first item, it is enough to show that whenever $(P,q)\in \Rc(T)$, the value $L_T(P,q)$ is a local minimum of $L_T(P,\cdot)$. Let $v_1,\dots, v_n$ be the terminals of $T$, let $\tilde{T}_1,\dots,\tilde{T}_l$ be the connected components of $T\smallsetminus\{v_1,\dots, v_n\}$ containing Steiner points and let $T_i$ be the closure of $\tilde{T}_i$ in $T$. Then each $T_i$ is a full topology (with terminals being a subset of terminals of $T$). Let $f: T\to \mathbb R^2$ be the embedding corresponding to $(P,q)$, then it is easy to see that the differential condition $\nabla_q L_T(p,q) = 0$ is equivalent to the fact that $f(T_i)$ is a locally minimal tree, which implies the statement. 

    For the second item, consider the mapping $M: \PP_2\to \PP_2\times (\mathbb R^2)^k$ which sends $P$ to the $(P,q(P))$. Then $M$ is continuous. But since $P_T(\Rc(T)) = M^{-1}(\Rg(T))$ and $\Rg(T)$ is open in $\PP_2\times (\mathbb R^2)^k$, we conclude that $\Rc(T)$ is open.
\end{proof}

Recall that  $P_T: \Rg(T)\to \PP_2$ is the projection. Given two topologies $T_1,T_2$ with $n$ labelled vertices, define
\[
    \Aa_{T_1,T_2} = \{ P\in P_{T_1}(\Rc(T_1))\cap P_{T_2}(\Rc(T_2))\ \mid\ L_{T_1}(P,q_1) = L_{T_2}(P,q_2), \text{ where } (P,q_i)\in \Rc(T_i) \}.
\]

In the next lemma we will use the notion of a subanalytic subset introduced in Section~\ref{subsec:Subanalytic subsets of a real analytic manifold}.

\begin{lemma}
    \label{lemma:AaT1T2_is_subanalytic}
    Let $T_1\neq T_2$ be two generic topologies with $n$ terminals. Then there exists an open set $U\subset \PP_2$ such that $\Aa_{T_1,T_2}$ is a subanalytic subset of $U$. In particular, $\Aa_{T_1,T_2}$ is a union of a countable collection of connected analytic submanifolds of $U$.
\end{lemma}
\begin{proof}
    Recall that for any $T$ the map $P_T:\Rg(T)\to \PP_2$ is the projection. Define $U = P_{T_1}(\Rc(T_1))\cap P_{T_2}(\Rc(T_2))$. By Lemma~\ref{lemma:properties_of_Rc} we have that $U$ is an open subset of $\PP_2$, and we have $\Aa_{T_1,T_2}\subset U$ by the definition of $\Aa_{T_1,T_2}$.
    
    Introduce the temporary notation $\Rc(T_i)_U = \Rc(T_i)\cap P_{T_i}^{-1}(U)$ and $\Rg(T_i)_U = \Rg(T_i)\cap P_{T_i}^{-1}(U)$ for simplicity. Recall the definition of a fiber product introduces in~\eqref{eq:def_of_fiber_product}. As we can see from the definition, 
    \[
        \Rg(T_1)_U\times_{P_{T_1} = P_{T_2}}\Rg(T_2)_U \subset U\times (\mathbb R^2)^{k_1}\times (\mathbb R^2)^{k_2}
    \]
    is and open subset, where $k_i$ is the number of Steiner points of $T_i$ and we identify $U\times (\mathbb R^2)^{k_1}\times (\mathbb R^2)^{k_2}$ with $(U\times (\mathbb R^2)^{k_1})\times_{P_{T_1} = P_{T_2}}(U\times (\mathbb R^2)^{k_2})$. Therefore $\Rg(T_1)_U\times_{P_{T_1} = P_{T_2}}\Rg(T_2)_U$ is a real analytic submanifold of $\Rg(T_1)_U\times \Rg(T_2)_U$. The set $\Rc(T_i)_U$ is an analytic subset of $\Rg(T_i)_U$, hence by Lemma~\ref{lemma:fiber_product} $\Rc(T_1)_U\times_{P_{T_1} = P_{T_2}}\Rc(T_2)_U$ is an analytic subset of $\Rg(T_1)_U\times_{P_{T_1} = P_{T_2}}\Rg(T_2)_U$.
    
    Define now
    \[
        \Rr = \{(P,q_1,q_2)\in \Rg(T_1)_U\times_{P_{T_1} = P_{T_2}}\Rg(T_2)_U\ \mid\ L_{T_1}(P,q_1) = L_{T_2}(P,q_2)\}.
    \]
    Then $\Rr$ is an analytic subset of $\Rg(T_1)_U\times_{P_{T_1} = P_{T_2}}\Rg(T_2)_U$. Denote by $\pi$ the natural projection $\pi: \Rg(T_1)_U\times_{P_{T_1} = P_{T_2}}\Rg(T_2)_U\to U$. Then we have
    \[
        \Aa_{T_1,T_2} = \pi\Bigl(\Rr\cap \Bigl(\Rc(T_1)_U\times_{P_{T_1} = P_{T_2}}\Rc(T_2)_U\Bigr)\Bigr).
    \]
    It follows from Lemma~\ref{lemma:images_of_semianalytic_are_subanalytic} that $\Aa_{T_1,T_2}$ is a subanalytic subset of $U$.

    The last assertion of the lemma follows from Corollary~\ref{cor:local_structure_of_subanalytic_sets}.
\end{proof}

\begin{lemma}
    \label{lemma:trees_in_IntAaT1T2_are_codirected}
    Let $T_1\neq T_2$ be two generic topologies with $n$ terminals, and assume that $P\in \Int (\Aa_{T_1,T_2})$ (here $\Int$ stays for the interior in $\PP_2$) and $S_1(P), S_2(P)$ are the images of $T_1,T_2$ on the plane. Then for each terminal $v$ we have $\deg_{T_1} v = \deg_{T_2} v$ and the trees $S_1(P), S_2(P)$ are codirected at $v$.
\end{lemma}
\begin{proof}
    Given $P\in \Int\Aa_{T_1,T_2}$, denote by $q_1(P)$ and $q_2(P)$ the configurations of Steiner points such that we have $(P,q_i(P))\in \Rc(T_i),\ i = 1,2$. Note that it is enough to prove the lemma for the points $P$ from a dense subset of the interior, since $q_i(P)$ depends continuously on $P$ (cf. the proof of Lemma~\ref{lemma:properties_of_Rc}). Recall that we denote by $P_{T_i}: \Rg(T_i)\to \PP_2$ the projection and $P_{T_i}(\Rc(T_i))$ is open in $\PP_2$; recall also that $P_{T_i}$ restricted to $\Rc(T_i)$ is one-to-one by Lemma~\ref{lemma:properties_of_Rc} and $q_i$ is the inverse mapping. From Lemma~\ref{lemma:local_structure_of_subanalytic_sets} and the fact that $P_{T_i}$ is one-to-one restricted to $\Rc(T_i)$ we find out that $q_i$ is differentiable on $P_{T_i}(\Rc(T_i))$ outside a subset which is nowhere dense in $P_{T_i}(\Rc(T_i))$. Thus, deforming $P$ inside $\Int (\Aa_{T_1,T_2})$ a little bit we can achieve that that both $q_1$ and $q_2$ are differentiable in a neighborhood of $P$.

    Further, we claim that deforming $P$ a little bit more we can assume that for any terminal vertex of degree 2 in $S_i(P)$ the angle between the corresponding edges in $S_i(P)$ is not equal to $\pi$ or $2\pi/3$. Indeed, let $v$ be such a vertex in, say, $S_1(P)$. Then $v$ divides $S_1(P)$ into two subtrees $S_1^+(P)$ and $S_1^-(P)$. Rotating $S_1^-(P)$ around $v$ a little bit we can assure that the angle at $v$ in $S_1(P)$ is not equal to $\pi$ or $2\pi/3$, and $q_1,q_2$ are still differentiable in a neighborhood of $P$. Repeating this for all terminal vertices of degree 2 in $S_1(P)$ and $S_2(P)$ we get the result; note that directions of edges in $S_1(P),S_2(P)$ depend on $P$ continuously, because $q_1(P),q_2(P)$ are continuous functions of $P$.

    Given a point $v$ from $P$ and an oriented edge $\vec{e}$ of $S_i(P)$ emanating from $v$ in $S_i(P)$ denote by $\eta_i(v, \vec{e})\in \mathbb R^2$ the unit vector in $\mathbb R^2$ codirected with $\vec{e}$. Set 
    \[
        \eta_i(v) = \sum_{\vec{e}\in \vec{E}(S_i(P))\colon\ o(\vec{e}) = v}\eta_i(v,\vec{e});
    \] note that the sum consists of one or two elements for each $v$ since all terminal vertices have degrees 1 or 2. A direct computation using angle condition at Steiner points of $S_i(P)$ shows that for any $\mu\in (\mathbb R^2)^n$ the derivative in the direction $\mu$ of $L_{T_i}(P, q_i(P))$ is given by
    \[
        \frac{\partial}{\partial \mu}L_{T_i}(P, q_i(P)) = -\sum_{j = 1}^n \eta_i(v_j)\cdot \mu_j,
    \]
    where $v_j$ is the $j$-th terminal. Since $L_{T_1}(\cdot,q_1(\cdot))$ and $L_{T_2}(\cdot, q_2(\cdot))$ are equal in a neighborhood of $P$, we conclude that $\eta_1(v_j) = \eta_2(v_j)$ for any $j = 1,\dots, n$. If $\deg_{S_1(P)} v_j = \deg_{S_2(P)}v_j = 1$, then this means that $S_1(P)$ and $S_2(P)$ are codirected at $v_j$. Assume that $\deg_{S_1(P)} v_j = 2$; then $|\eta_1(v_j)|\neq 1$ since the angle between the two edges emanating from $v_j$ is not equal to $2\pi/3$ by our assumption. From this inequality and the equality $\eta_1(v_j) = \eta_2(v_j)$ we find out that $\deg_{S_2(P)} v_j = 2$ also. Further, we have $\eta_1(v_j)\neq 0$ since the angle between the two edges emanating from $v_j$ is not equal to $\pi$ by our assumption. As a consequence, there is only one unordered pair of unit vectors $(\mu_+,\mu_-)$ such that $\eta_1(v_j) = \mu_++\mu_-$, hence the pair of edges emanating from $v_j$ in both $S_1(P)$ and $S_2(P)$ must have these directions, so that $S_1(P)$ and $S_2(P)$ are codirected at $v_j$. We conclude that $S_1(P)$ and $S_2(P)$ are codirected at all terminals.
\end{proof}

Let us now formulate the following theorem of Oblakov:
\begin{theo}[Oblakov~\cite{oblakov2009non}]
    \label{theo:oblakov}
    Assume that $S_1$ and $S_2$ are two locally minimal trees connecting the same set of terminals $P\in \PP_2$ and codirected at this set. Then $S_1$ and $S_2$ coincide.
\end{theo}

Given two generic topologies $T_1$ and $T_2$ with $n$ terminals and $P\in \Aa_{T_1,T_2}$, denote by $S_1(P)$ and $S_2(P)$ the embeddings of $T_1$ and $T_2$ as in the lemma above. Define
\[
    \Aa_{T_1,T_2}^{\mathrm{min}} = \{P\in \Aa_{T_1,T_2} \ \mid\ \text{$S_1(P)$ and $S_2(P)$ are both locally minimal trees}\}.
\]
We have the following 

\begin{cor}
    \label{cor:interior_of_AaT1T2_not_intersect_Aa}
    Let $T_1\neq T_2$ be two generic topologies with $n$ terminals. Then we have
    \[
        \Int (\Aa_{T_1,T_2})\cap \Aa_{T_1,T_2}^{\mathrm{min}} = \varnothing,
    \]
    where $\Int$ stays for the interior in $\PP_2$.
\end{cor}
\begin{proof}
    Follows from Lemma~\ref{lemma:trees_in_IntAaT1T2_are_codirected} and Theorem~\ref{theo:oblakov}.
\end{proof}

\begin{lemma}
    \label{lemma:Aa_in_another}
    We have
    \[
        \Aa\subset \bigcup_{\substack{T_1\neq T_2\\ T_1,T_2\text{ - generic}}} \Aa^{\mathrm{min}}_{T_1,T_2}\cup \{ P\in \PP_2\ \mid\ \text{there is a Steiner tree for $P$ with a non-generic topology} \}.
    \]
\end{lemma}
\begin{proof}
    Let $P\in \Aa$, then $P$ is connected by two Steiner trees $S_1$ and $S_2$. If both $S_1$ and $S_2$ have the same topology $T$, then we get the contradiction with Lemma~\ref{lemma:properties_of_Rc}. Thus the topologies of $S_1$ and $S_2$ are different and the lemma follows.
\end{proof}

We can now proof Theorem~\ref{main}:

\begin{proof}[Proof of Theorem~\ref{main}]
    By Lemma~\ref{lemma:dim_of_Amb} the dimension of the set $\Aa$ of ambiguous configurations is at least $2n-1$. Therefore, by Lemma~\ref{lemma:Aa_in_another} and Observation~\ref{obs:topologies},~(iv), we have
    \[
        \dim \Aa \leq \max\{\dim \Aa^{\mathrm{min}}_{T_1,T_2} \ \mid\ T_1\neq T_2,\ T_1,T_2\text{ - generic}\}.
    \]
    Let two generic topologies $T_1\neq T_2$ be fixed. Let $P\in \Aa_{T_1,T_2}^{\mathrm{min}}$. By Lemma~\ref{lemma:AaT1T2_is_subanalytic} we have $\Aa_{T_1,T_2} = X_1\cup X_2\cup\ldots$, where $X_1,X_2,\dots$ are connected analytic submanifolds of an open subset $U\subset \PP_2$. Therefore, by Corollary~\ref{cor:interior_of_AaT1T2_not_intersect_Aa}
    \[
        \Aa^{\mathrm{min}}_{T_1,T_2}\subset \bigcup_{i\ \colon\ \dim X_i \leq 2n-1}X_i.
    \]
    It follows that $\dim \Aa^{\mathrm{min}}_{T_1,T_2}\leq 2n-1$. Since $T_1,T_2$ were arbitrary, we conclude that $\dim \Aa\leq 2n-1$.
\end{proof}

\begin{remark}
    A reasonable question would be if the dimension of the whole $\Aa_{T_1,T_2}$ (not only $\Aa_{T_1,T_2}^{\min}$) is less or equal to $2n-1$. We claim that it can be proven in the similar way we prove it for $\Aa_{T_1,T_2}^{\min}$. Indeed, due to Lemma~\ref{lemma:AaT1T2_is_subanalytic} it is enough to prove that the interior of $\Aa_{T_1,T_2}$ is empty. If $P\in \Int\Aa_{T_1,T_2}$ and $S_1(P)$ and $S_2(P)$ are the corresponding embeddings of $T_1,T_2$, then by Lemma~\ref{lemma:trees_in_IntAaT1T2_are_codirected} trees $S_1(P)$ and $S_2(P)$ are codirected. Our claim is now that Theorem~\ref{theo:oblakov} can still be applied in this case to say that $S_1(P)$ and $S_2(P)$ coincide. Note that $S_1(P)$ and $S_2(P)$ are non necessary locally minimal networks as we allow them to have angles less than $\frac{2\pi}{3}$ at terminals of degree 2. Nevertheless, the proof of Theorem~\ref{theo:oblakov} given by Oblakov~\cite{oblakov2009non} still applies to this case.
\end{remark}

\section{Steiner trees in real analytic Riemannian manifolds}
\label{sec:Steiner trees in analytic Riemannian manifolds}

A question on the uniqueness of Steiner trees in a Riemannian manifold was raised in~\cite{edelsbrunner2012current}.

We should not expect that Theorem~\ref{main} can be directly generalized to the case when $\mathbb R^2$ is replaced with an arbitrary manifold $M$ (cf. Section~\ref{subsec:topological_example}). Nevertheless, if $M$ is a real analytic manifold, then we still can expect that the set of ambiguous configurations of terminals either has a non-empty interior, or dimension strictly less than the set of all configurations of terminals.

The aim of this section is to build a similar framework to the one used in the proof of Theorem~\ref{main} in Section~\ref{subsec:proof_of_Theorem_on_dimension} in the case of arbitrary real analytic manifold $M$. Using this framework we reduce the alternative stated above to Conjecture~\ref{conjecture:Rm_is_semianalytical} about analytic sets.

\subsection{Realization space $\Rg$ for an arbitrary metric space}
\label{subsec:arbitrary_metric_space}

Let us begin with rephrasing Problem~\ref{Problem1} with an arbitrary proper metric space $M$ instead of $\mathbb R^d$.

\begin{problem}\label{problem2}
Let $M$ be a metric space.
For a given finite set $Q = \{p_1,\dots ,p_n\} \subset M$ find a connected set $\St$ with minimal length (one-dimensional Hausdorff measure) containing $Q$.
\end{problem}

Due to the following theorem solutions to Problem~\ref{problem2} still lie among geodesically embedded trees:

\begin{theo}[Paolini--Stepanov,~\cite{paolini2013existence}]
Assume that $M$ is proper (i.e. all closed balls in $M$ are compact) and connected. Then a solution of Problem~\ref{problem2} exists.
Moreover, for any solution $\St(Q)$ the following statements hold:
\begin{itemize}
    \item [(i)] $\St$ is compact;
    \item [(ii)] $\St$ contains no closed loops (homeomorphic images of $\mathbb S ^ 1$);
    \item [(iii)] $\St \setminus Q$ has a finite number of connected components, each component has strictly positive length, and the closure of each component is a finite geodesic embedded graph with endpoints in $Q$;
    \item [(iv)] the closure of every connected component of $\St \setminus Q$ is a topological tree with
endpoints in $Q$, and all the branching points having finite degree.
\end{itemize}
\end{theo}

This theorem motivates the following definition. Given a positive integer $n$ and a topology $T$, define
\begin{equation}
    \Rg(T) = \{ f: T\to M\ \mid\ f\text{ is an embedding which maps all edges to (shortest) geodesics} \}/_\sim,
\end{equation}
where $f_1\sim f_2$ if and only if $f_1(v) = f_2(v)$ for any labelled vertex $v$ and $f_1(T) = f_2(T)$ as subsets of $M$. Note that this $\Rg(T)$ is a straightforward generalization of $\Rg(T)$ introduces in Section~\ref{subsec:proof_of_Theorem_on_dimension}. Let $\mathcal{K}(M)$ be the set of all compact subsets of $M$ endowed with the Hausdorff metric. Introduce two maps: first, $\Rg(T)\to \mathcal{K}(M)$ which maps $f$ to $f(T)$, second, $P_T: \Rg(T)\to M^n$ which sends $f$ to the collection $(f(v_1),\dots,f(v_n))$ of images of $n$ terminals. The topology on $\Rg(T)$ is defined to be the pullback of the product topology under the map $\Rg(T)\to \mathcal{K}(M)\times M^n$ given by the two maps above (the map $P_T$ is added to keep track of the enumeration of terminals).

\subsection{Manifold structure on \texorpdfstring{$\Rg$}{Rgeomb}}

Let us now assume that $M$ is a connected real analytic manifold with a Riemannian metric $d$ which depends analytically on the point of $M$. We define the intrinsic metric $d_{in}$ on $M$ as usual; note that $(M, d_{in})$ is a proper metric space. Given a point $p$, denote by $\exp_p$ the exponential map defined with respect to $d$; since we have not required $(M, d_{in})$ to be complete, $\exp_p$ is defined only for an open subset of the tangent space $T_pM$. Set
\[
    \widetilde{T_pM} = \{ w\in T_pM\ \mid\ w\neq 0\text{ and $\exp_p(w)$ is defined} \}.
\]
Denote by $TM$ the tangent bundle of the manifold $M$. Elements of $TM$ are parameterized by pairs $(p, w)$ where $p\in M$ and $w\in T_pM$. Let $\widetilde{TM}\subset TM$ be the union of $\widetilde{T_pM}$ over all $p\in M$. Then $\widetilde{TM}$ is an open subset of $TM$ and $\exp: \widetilde{TM}\to TM$, given by $(p,w)\mapsto (\exp_p(w), \mathrm{d}\exp_p(w))$ on each fiber, is a real analytic map mapping $\widetilde{TM}$ onto its image in $TM$ diffeomorphically (see~\cite[Section 8]{coddington1955theory}). Moreover, this diffeomorphism is analytic at any point.

Let us show that for any topology $T$ the set $\Rg(T)$ has a natural structure of an analytic manifold. Let $e_1,\dots,e_m$ be the edges of $T$ and $\vartheta$ be an arbitrary orientation on edges of $T$. Define the map $\varphi_{\vartheta}: \Rg(T)\to (\widetilde{TM})^m$ by
\[
    \varphi_{\vartheta}(f) = ((f(o(e_1)), w_1),\dots, (f(o(e_m)), w_m))\in (\widetilde{TM})^m,
\]
where $o(e_k)$ is the origin of $e_k$ oriented according to $\vartheta$ and $w_k$ is the tangent vector at $o(e_k)$ such that the geodesic $f(e_k)$ is given by $\{\exp_{o(e_k)}(tw_k)\}_{0\leq t\leq 1}$.

\begin{lemma}
    \label{lemma:coordinates_on_E}
    The following statements hold true:
    \begin{itemize}
        \item[(i)] For any orientation $\vartheta$, the map $\varphi_{\vartheta}$ is a homeomorphism between $\Rg(T)$ and a smooth real analytic submanifold $\mathcal{U}_{\vartheta}\subset (\widetilde{TM})^m$ of dimension $(m+1)\cdot \dim M$.
        \item[(ii)] The morphism $P_T\circ \varphi_{\vartheta}^{-1}: \mathcal{U}_{\vartheta}\to M^n$ is analytic.
        \item[(iii)] Given two orientations $\vartheta_1, \vartheta_2$, the map $\varphi_{\vartheta_2}\circ\varphi_{\vartheta_1}^{-1}: \mathcal{U}_{\vartheta_1}\to \mathcal{U}_{\vartheta_2}$ is analytic.
    \end{itemize}
\end{lemma}
\begin{proof}
    We first prove (i). Note that $\varphi_\vartheta$ is injective. Let us show now that $\mathcal{U}_\vartheta = f(\Rg(T))$ is a real analytic submanifold. In fact we have
    \begin{equation}
    \label{eq:mathcalU}
        \mathcal{U}_\vartheta = \{ ((x_1,w_1),\dots,(x_m,w_m))\in (\widetilde{TM})^m\ \mid\ \exp_{x_i}(w_i) = x_j\text{ if }x_j = \mathrm{tail}(e_i),\ x_i = x_j\text{ if }o(e_i) = o(e_j) \},
    \end{equation}
    where $\mathrm{tail}(e_i)$ is the tail of $e_i$ oriented according to $\vartheta$. It follows that $\mathcal{U}_\vartheta$ is an analytic subset of $(\widetilde{TM})^m$. The smoothness easily follows from the fact that $\exp$ is a diffeomorphism.

    Finally, we have that the inverse mapping $\varphi_\vartheta^{-1}:\mathcal{U}_\vartheta\to \Rg(T)$ is continuous because $\exp$ is continuous. We conclude that $\varphi_\vartheta$ is continuous because $\mathcal{U}_\vartheta$ is locally compact.

    The items (ii) and (iii) follows easily from that fact that $\exp$ is analytic.
\end{proof}

From Lemma~\ref{lemma:coordinates_on_E} we see that $\Rg(T)$ has a natural structure of a real analytic manifold.
\begin{lemma}
    \label{lemma:projection_P_analytic}
    The map $P_T: \Rg(T)\to M^n$ is analytic and $P_T(\Rg(T))$ is open in $M^n$. The differential $\mathrm{d}P$ has the maximal rank at any point. In particular, the fiber of $P$ is smooth and has the dimension $\dim \Rg(T) - n\dim M$. 
\end{lemma}
\begin{proof}
    Choose an arbitrary orientation $\vartheta$ on the edges of $T$ having the property that any terminal is an origin of some edge. The lemma now follows from Lemma~\ref{lemma:coordinates_on_E} and the description~\eqref{eq:mathcalU} of $\mathcal{U}_\vartheta$.
\end{proof}

\subsection{The subvariety \texorpdfstring{$\Rm$}{Rlocmin} of \texorpdfstring{$\Rg$}{Rgeomb}}
\label{subsec:Rlocmin}

Define the length function $L_T: \Rg(T)\to \mathbb R$ by
\[
    L_T(f) = \mathcal{H}^1(f(T)).
\]
Let $\Rm(T)\subset \Rg(T)$ be the subset given by
\[
    \Rm(T) = \{ f\in \Rg(T)\ \mid\ f \text{ is a local minimum of $L_T$ on the set $P_T^{-1}(P_T(f))$} \}.
\]

Note that, given the orientation $\vartheta$ on the edges of $T$ we have
\[
    L\circ \varphi_{\vartheta}^{-1}((x_1,w_1), \dots, (x_m,w_m)) = |w_1| + \dots + |w_m|
\]
(cf. Lemma~\ref{lemma:coordinates_on_E}), hence $L_T$ is an analytic function on $\Rg(T)$. Recall that due to Lemma~\ref{lemma:projection_P_analytic}, the fibers of $P_T$ are smooth. Define the vertical gradient of the function $L_T$ at a point $f\in \Rg$ to be the restriction of $\nabla L_T$ to the tangent space to the fiber $P_T^{-1}(P_T(f))$ of $P_T$ over $f$. Define
\[
    \Rc(T) = \{ f\in \Rg\ \mid\ \text{vertical gradient of $L_T$ at $f$ is zero} \}.
\]
Clearly, $\Rm(T)\subset \Rc(T)$; note that if $M = \mathbb{R}^2$ with the Euclidean metric, then the converse inclusion is also true, but in general we do not have equality between these two sets. We expect nevertheless that $\Rm(T)$ is a subanalytic subset of $\Rg(T)$. By Lemma~\ref{lemma:projection_P_analytic}, this statement would immediately follow from the following assertion:

\begin{conjecture}
    \label{conjecture:Rm_is_semianalytical}
    Assume that $U\subset \mathbb R^m$ and $V\subset \mathbb R^k$ are open subsets, $f: U\times V\to \mathbb R$ is real analytic. Put
    \[
        A = \{ (x,y)\in U\times V\ \mid\ x\text{ is a local minimum for $f(\cdot, y)$ on $U$, if $y$ is fixed} \}.
    \]
    Then $A$ is a semianalytic subset of $U\times V$.
\end{conjecture}

Assume that for any $T$ the set $\Rm(T)$ is subanalytic for a moment. In this case we immediately have the following

\begin{prop}
    \label{prop:ambiguous_crit_conf}
    Let $T_1,T_2$ be two topologies with $n$ terminals and assume that  
    \[
        \Aa_{T_1,T_2} = \{ P\in M^n\ \mid\ \exists f_1\in\Rm(T_1),\ f_2\in \Rm(T_2)\ \colon\ P_{T_1}(f_1) = P_{T_2}(f_2) = P\text{ and }L_{T_1}(f_1) = L_{T_2}(f_2) \},
    \]
    Then either $\Aa_{T_1,T_2}$ has non-empty interior, or the Hausdorff dimension of $\Aa_{T_1,T_2}$ is strictly less than $n\dim M$.
\end{prop}
Note that both alternatives in Proposition~\ref{prop:ambiguous_crit_conf} can occure, see Section~\ref{subsec:topological_example}.
\begin{proof}
    The proof essentially follows proof of Lemma~\ref{lemma:AaT1T2_is_subanalytic}. We consider the fiber product $\Rm(T_1)\times_{P_{T_1} = P_{T_2}}\Rm(T_2)$ which is subanalytic due to Lemma~\ref{lemma:fiber_product}. Inside $\Rm(T_1)\times_P\Rm(T_2)$ we consider the set $\tilde{\Aa}_{T_1,T_2}$ cut out by the equation $L_{T_1}(f_1) = L_{T_2}(f_2)$. Then $\tilde{\Aa}_{T_1,T_2}$ is subanalytic and $\Aa_{T_1,T_2} = P(\tilde{\Aa}_{T_1,T_2})$, where $P = P_{T_1}: \Rm(T_1)\times_{P_{T_1} = P_{T_2}}\Rm(T_2)\to M^n$ is the projection. The proposition now follows from Corollary~\ref{cor:local_structure_of_subanalytic_sets}.
\end{proof}

Let us finalize this section with a short discussion on Conjecture~\ref{conjecture:Rm_is_semianalytical}. 
Note that if $k = 0$ so that $V$ is just a point, then the statement in the conjecture is straightforward (see also~\cite{FernandoLocMin} for much more general statement). Indeed, define first
\[
    Z = \{ x\in U\ \mid\ \nabla f(x) = 0 \},
\]
Clearly, $A\subset Z$. Using~\cite[Proposition~2.10]{bierstone1988semianalytic} we find that each $(x,y)\in U$ has a neighborhood $U'$ such that $Z\cap U' = Z_1\cup\dots\cup Z_k$ where each $Z_j$ is a connected semianalytic subset of $U'$ which is also an analytic submanifold of $U$. In particular, $f$ restricted to each $Z_j$ is a constant, say, $C_j$. Define
\[
    B_j = \overline{\{ x\in U\ \mid\ f(x) < C_j \}}.
\]
Then $B_j$ is a semianalytic subset (as the closure of a semianalytic subset is semianalytic~\cite[Proposition~2.10]{bierstone1988semianalytic}).
We have
\[
    A\cap U' = \bigcup_{j = 1}^k \left( Z_j\smallsetminus B_j \right).
\]
Since $x$ was chosen arbitrary, this proves that $A$ is semianalytic.

With some more involved arguments using Weierstrass preparation theorem we can prove the conjecture when $\dim U = 1$ and $V$ is arbitrary, but the general case remains unclear to us.

\subsection{Example of \texorpdfstring{$\Aa_{T_1,T_2}$}{AT1T2} with non-empty interior}
\label{subsec:topological_example}

In this subsection we construct an example of a Riemann surface $M$ with a locally flat metric such that there exist two topologies $T_1,T_2$ on 8 terminals for which $\Aa_{T_1,T_2}$ has non-empty interior (see Proposition~\ref{prop:ambiguous_crit_conf} for the definition of $\Aa_{T_1,T_2}$).

Let $T_1$ and $T_2$ be the two topologies introduced by Ivanov and Tuzhilin in~\cite[Fig.~1]{ivanov2006uniqueness}, see Figure~\ref{4exex}. Let $T_1$ be the topology of the tree drawn by solid lines for certainty. We fix a set of terminals $x = (x_1,\dots,x_8)\in \PP_2$ and fix \emph{immersions} $f_i: T_i\to \mathbb{R}^2$ into plane as on Figure~\ref{4exex}. Following~\cite{ivanov2006uniqueness} we have the following

\begin{lemma}
    \label{lemma:non-uniqueness}
    For any configuration $(\tilde{x}_1,\dots,\tilde{x}_8)$ sufficiently close to $(x_1,\dots,x_8)$ the corresponding immersions $\tilde{f}_1,\tilde{f}_2$ of the topologies $T_1,T_2$ realizing them as locally minimal trees are codirected and have the same length.
\end{lemma}
\begin{proof}
    Follows immediately from Melzak algorithm. Note that being codirected at some point implies being codirected in a neighborhood, which in turn is equivalent to the equality of length (cf. Lemma~\ref{lemma:trees_in_IntAaT1T2_are_codirected}).
\end{proof}

Let 
\[
    T = T_1\sqcup T_2/_{\text{glued along terminal half-edges}},
\]
i.e. we glue $T_1$ and $T_2$ along the portion of edges emanating from $x_1,\dots, x_8$ which coincide on the picture (other intersection points on the picture are not glued). Let $f: T\to \mathbb R^2$ be the correspondig immersion.

Note that $T$ inherits a metric from $f(T)$: use Euclidean metric on $f(T)$ to measure distances along an edge of $T$, and use the inner metric to measure the distance between two arbitrary points.

\begin{figure}[h]
    \centering
    \begin{tikzpicture}[scale=1.5]

    \begin{scope}[shift={(0,0)}]
        \draw[ultra thick, blue] (0,0) coordinate (x1) --++  (-60:1) coordinate (B) --++ (0:2) coordinate (C) --++ (60:1) coordinate (x2);
        \draw[ultra thick, blue] (C) --++ (-60:1.5) coordinate (D) --++ (0:1) coordinate (x3);
        \draw[ultra thick, blue] (D) --++ (-120:1.5) coordinate (E) --++ (-60:1.5) coordinate (x5);
        \draw[ultra thick, blue] (E) --++ (180:3) coordinate (x7);
         \draw[ultra thick, blue] (B) --++ (-120:1) coordinate (A) --++ (180:1) coordinate (x8);
        \draw[ultra thick, blue] (A) --++ (-60:2.5) coordinate (F) --++ (-120:1) coordinate (x6);
         \draw[ultra thick, blue] (F) --++ (0:3) coordinate (x4);

        \draw[ultra thick, dashed, black] (x8) --++ (0:3) coordinate (a) --++ (60:2);
       \draw[ultra thick, dashed, black] (a) --++ (-60:2.5) coordinate (b) --++ (0:1);
       \draw[ultra thick, dashed, black] (b) --++ (-120:0.5) coordinate (c) --++ (-60:0.5);
       \draw[ultra thick, dashed, black] (c) --++ (180:2) coordinate (d) --++ (-120:0.5);
       \draw[ultra thick, dashed, black] (d) --++ (120:1) coordinate (e) --++ (180:1);
       \draw[ultra thick, dashed, black] (e) --++ (60:1.5) coordinate (f) --++ (0:3);
       \draw[ultra thick, dashed, black] (f) --++ (120:2.5);

       \foreach \x in {1,2,3,4,5,7,8}{
             \fill (x\x) circle (2pt) node[above right] {$x_\x$};
        }
        \fill (x6) circle (2pt) node[above left] {$x_6$};
       
  %    \fill (A) circle (2pt) node[right]{$A$};
  %      \fill (B) circle (2pt) node[right]{$B$};
  %      \fill (C) circle (2pt) node[right]{$C$};
    \end{scope}

\end{tikzpicture}
    \caption{Two locally minimal trees with self-intersections}
    \label{4exex}
\end{figure}

Fix an $\eps>0$ and for each $t\in T$ define the surface $M_t$ to be a copy the $\eps$-neighborhood of $f(t)$ in $\mathbb{R}^2$. The map $f$ extends to the embedding $F_t:M_t\to \mathbb R^2$. Given $t_1,t_2\in T$ on the distance at most $10\eps$ from each other, glue $M_{t_1}$ and $M_{t_2}$ such that $F_{t_1} = F_{t_2}$. As a result we obtain Riemann surface $M$ containing $T$, and an immersion $F: M\to \mathbb{R}^2$ such that the image $F(M)$ is the $\eps$-neighborhood of $f(T)$. Note that $M$ is endowed with a locally flat metric for which $F$ is a local isometry. 

Let $g_i: T_i\to M, i = 1,2,$ be the natural maps. Note that $f_i = F\circ g_i$ and $g_i$ are injective. Denote by $X_1,\dots, X_8\in M$ the points such that $F(X_i) = x_i$. The following proposition is straightforward.

\begin{prop}
    \label{prop:non-uniqueness}
    Trees $g_1(T_1), g_2(T_2)$ are locally minimal on $M$ and $(g_1,g_2)\in \Int \Aa_{T_1,T_2}$.
\end{prop}

We are not able to extend this example on minimal trees.

\bigskip

\textbf{Acknowledgements.} The work is supported by Russian Science Foundation grant 14-21-00035. The authors thank Roman Karasev for pointing out the reference~\cite{bierstone1988semianalytic} that helped us to simplify our proofs significantly and make them more consistent.

\bibliographystyle{plain}
\bibliography{main}

\section*{Appendix}

\begin{proof}[Proof of Lemma~\ref{lemma:class_of_embeddings}]
    We construct a map $F_i: PM_i\to PM_{i+1}$ for $i = 1,2,3$ (where we compute indices$\mod 3$), and then show that the composition of these maps is identity.

    Let $(T,\sigma)\in PM_1$ be given, let $N = |\vec{E}(T)|$. If $N=0$, then we just take an empty $D$. Assume that $N>0$. Define $\alpha: \vec{E}(T)\to \vec{E}(T)$ to be the involution reversing the orientation and set $\varphi = \alpha\circ \sigma$. It is easy to verify that $\varphi$ is a cyclic permutation of $\vec{E}(T)$. Set $f:\vec{E}(T)\to \partial D$ to be any bijection which respects the cyclic order imposed by $\varphi$. It is easy to check that $(T,[f])$ belongs to $PM_2$, thus we get the map $F_1$.

    Let $(T,[f])\in PM_2$ be given. If $T$ is one point, then we define $\iopta$ arbitrary. Assume that $T$ has at least two vertices, let $N = |\vec{E}(T)|$ and $D$ be the regular $N$-gon. As above, let $\alpha: \vec{E}(T)\to \vec{E}(T)$ to be the involution reversing the orientation. Now, glue each edge $\vec{e}\in \partial D$ with $f\circ \alpha\circ f^{-1}(\vec{e})$ in opposite direction. It is straightforward to see that in this way we get an oriented surface $S$ out of $D$, and a natural embedding $\iopta$ of $T$ into $S$. Computing the Euler characteristic we find out that $S$ is a sphere and hence $\iopta$ corresponds to a planar embedding of $T$. Set $F_2(T,[f]) = (T,[\iopta])$ (where the topology on $T$ comes from $\iopta$ naturally).

    Let $(T,[\iopta])\in PM_3$ be given. If $T$ is one point, then we take $\sigma$ to be the only map between empty sets. Assume that $T$ has at least two vertices. Let $v$ be a vertex and $\vec{E}_v(T)$ be the set of oriented edges emanating from $v$. Then given $\vec{e}\in \vec{E}_v(T)$ define $\sigma(\vec{e})$ to be first edge in $\vec{E}_v(T)$ coming after $\vec{e}$ when going around $\iopta(v)$ in counterclockwise direction. Set $F_3(T,[\iopta]) = (T,\sigma)$.

    The fact that $F_3\circ F_2\circ F_1 = \mathrm{id}$ is a simple exercise which we leave to the reader. Note that given a labelled tree $T$ the amount of all possible $\sigma$ such that $(T,\sigma)\in PM_1$, or $f$ such that $(T,[f])\in PM_2$ is finite; this shows that $F_1$ is a bijection. On the other hand, the fact that the number of homotopy classes of embeddings $\iopta$ for a given topology is finite is not obvious. Hence, at the moment we have only a right inverse for $F_3$. Let us sketch the construction the inverse map for $F_2$ to overcome this difficulty. Choose an embedding $\iopta:T\to \mathbb C$ and consider the simply-connected surface $\widehat{\mathbb{C}}\smallsetminus \iopta(T)$, where $\widehat{\mathbb{C}}$ is the Riemann sphere. Let $\mathbb D$ be the unit disc and $\psi: \mathbb D\to \widehat{\mathbb{C}}\smallsetminus \iopta(T)$ be the uniformization map. One can show that $\psi$ extends to the boundary of $\mathbb D$ in a unique way such that $\psi: \overline{\mathbb D}\to \widehat{\mathbb{C}}$ is continuous. Moreover, each point of $\iopta(T)$ corresponds to several prime ends of the domain $\widehat{\mathbb{C}}\smallsetminus \iopta(T)$ (see \cite[Chapter~2]{Pommerenke}); there are two prime ends for each inner point of an edge, and $\deg v$ prime ends for each vertex $v$. Let $v_1,\dots,v_N$ be all the preimages of vertices of $T$ on $\partial \mathbb D$, the count of the prime ends implies that $N = |\vec{E}(T)|$. Then $\overline{\mathbb D}$ together with these points has the combinatorics of the regular $N$-gon, whence we get the morphism $f$ such that $(T,f)\in PM_2$. The fact that this construction inverses $F_2$ is straightforward.
\end{proof}

\begin{proof}[Proof of Lemma~\ref{lemma:def_of_embedding_type}]
    Let $(T,\sigma)$ be given and $N = |\vec{E}(T)|$. Then the length of the word $C(T,\sigma)$ is $N$, hence $C$ distinguishes pairs $(T,\sigma)$ with different cardinality $N$ of the set of edges of the tree. We will show that $C$ distinguishes different pairs with the same $N$ by induction. If $N=0,2,4$, then there is nothing to prove, assume that $N>4$ and $W = C(T,\sigma)$. We need to show that if $W = C(T_1,\sigma_1)$, then $(T, \sigma) = (T_1,\sigma_1)$. Define 
    \begin{align*}
        &I_1 = \{i\ \mid\ a_i\text{ occurs 1 time in }W\}\\
        &I_2 = \{i\ \mid\ a_i\text{ occurs 2 times in }W\}.
    \end{align*}
    We clearly have a bijection between the labels $\{a_i\ \mid\ i\in I_1\}$ and $\{a_i\ \mid\ i\in I_2\}$ and the vertices of degree 1 and 2 in $T$ respectively, and the same for $T_1$. Assume that we can find $i\in I_1$ such that $i+1\in I_2$ (here $N+1=1$). Then consider the word $W'$ obtained from $W$ by removing $a_i$ and $a_{i+1}$. Then $W' = C(T',\sigma')$, where $T'$ is obtained from $T$ by removing the edge $a_ia_{i+1}$ and keeping all labels, and $\sigma'$ is computed from $\sigma$ in the natural way (note that $T$ has at least one edge since we assume that $N>4$). In the same time, $W' = C(T'_1,\sigma'_1)$, where $(T'_1,\sigma'_1)$ is obtained from $(T_1,\sigma_1)$ in the same procedure. By the induction hypothesis $(T',\sigma') = (T_1',\sigma_1')$. From here, it is easy to see that $(T,\sigma) = (T_1,\sigma_1)$.

    Assume now that for any $i\in I_1$ we have $i+1\notin I_2$. It follows that one can find $i\in I_1$ such that $i+2\in I_1$ also. Consider the word $W'$ obtained from $W$ by removing $a_{i-1},a_i,a_{i+1}$ and $a_{i+3}$. This word corresponds to $C(T',\sigma')$, where $T'$ is obtained from $T$ by removing two vertices $a_i$ and $a_{i+2}$ and labelling their common parent by $a_{i+2}$ (note that their parent must have degree 3). Note that $T'$ has at least one edge since $N>4$. Doing the same with $T_1$ we again get two pairs $(T',\sigma')$ and $(T'_1,\sigma'_1)$ such that $W' = C(T',\sigma') = C(T_1',\sigma_1')$, which implies that $(T',\sigma') = (T_1',\sigma_1')$ by the induction and, eventually, $(T,\sigma) = (T',\sigma')$.
\end{proof}

\end{document}